\documentclass[12pt,oneside,reqno]{amsart}
\usepackage{txfonts}
\usepackage{bbm}
\usepackage{amsmath}
\usepackage{graphicx}
\usepackage{mathrsfs}
\usepackage{stmaryrd}
\usepackage{amsfonts}
\usepackage{enumerate,amsmath,amssymb,amsthm}

\numberwithin{equation}{section}

\newcommand{\be}{\begin{eqnarray}}
\newcommand{\ee}{\end{eqnarray}}
\newcommand{\ce}{\begin{eqnarray*}}
\newcommand{\de}{\end{eqnarray*}}
\newtheorem{theorem}{Theorem}[section]
\newtheorem{lemma}[theorem]{Lemma}
\newtheorem{remark}[theorem]{Remark}
\newtheorem{definition}[theorem]{Definition}
\newtheorem{proposition}[theorem]{Proposition}

\newtheorem{corollary}[theorem]{Corollary}

\def\a{\alpha}
\def\om{\omega}

\def\p{\partial}

\def\eps{\epsilon}

\def\l{\lambda}

\def\[{{\Big[}}
\def\]{{\Big]}}
\def\<{{\langle}}
\def\>{{\rangle}}
\def\({{\Big(}}
\def\){{\Big)}}

\def\dif{{\mathord{{\rm d}}}}

\def\div{{\mathord{{\rm div}}}}

\def\u{\mathord{{\bf u}}}
\def\f{\mathord{{\bf f}}}
\def\v{\mathord{{\bf v}}}
\def\e{\mathord{{\bf e}}}
\def\je{\mathord{{\bf j}}_\epsilon}

\def\w{\mathord{{\bf w}}}
\def\h{\mathord{{\bf h}}}

\def\no{\nonumber}
\def\bt{\begin{theorem}}
\def\et{\end{theorem}}
\def\bl{\begin{lemma}}
\def\el{\end{lemma}}
\def\br{\begin{remark}}
\def\er{\end{remark}}
\def\bd{\begin{definition}}
\def\ed{\end{definition}}
\def\bp{\begin{proposition}}
\def\ep{\end{proposition}}
\def\bc{\begin{corollary}}
\def\ec{\end{corollary}}
\def\cA{{\mathcal A}}
\def\cB{{\mathcal B}}

\def\cE{{\mathcal E}}
\def\cF{{\mathcal F}}

\def\cJ{{\mathcal J}}
\def\cK{{\mathcal K}}

\def\cN{{\mathcal N}}
\def\cO{{\mathcal O}}
\def\cP{{\mathcal P}}

\def\cS{{\mathcal S}}

\def\cV{{\mathcal V}}

\def\mB{{\mathbb B}}

\def\mD{{\mathbb D}}
\def\mE{{\mathbb E}}

\def\mH{{\mathbb H}}

\def\mK{{\mathbb K}}

\def\mN{{\mathbb N}}

\def\mP{{\mathbb P}}

\def\mR{{\mathbb R}}

\def\mT{{\mathbb T}}
\def\mU{{\mathbb U}}

\def\mW{{\mathbb W}}
\def\mX{{\mathbb X}}

\def\geq{\geqslant}
\def\leq{\leqslant}

\def\sE{{\mathscr E}}

\def\sH{{\mathscr H}}

\def\sP{{\mathscr P}}

\def\bT{{\mathbf T}}

\def\lb{{\llbracket}}
\def\rb{{\rrbracket}}

\allowdisplaybreaks

\begin{document}

\title{Stochastic Tamed 3D Navier-Stokes Equations: Existence, Uniqueness and Ergodicity}

\date{}
\author{Michael R\"ockner$^{1,2}$, Xicheng Zhang$^{1,3,4}$ }

\dedicatory{
$^1$Fakult\"at f\"ur Mathematik,
Universit\"at Bielefeld\\
Postfach 100131,
D-33501 Bielefeld, Germany\\
$^2$Departments of Mathematics and Statistics, Purdue University\\
W. Laffayette, IN, 47907, USA\\
$^3$School of Mathematics and Statistics\\
The University of New South Wales, Sydney, 2052, Australia\\
$^4$Department of Mathematics,
Huazhong University of Science and Technology\\
Wuhan, Hubei 430074, P.R.China\\
 M. R\"ockner: roeckner@math.uni-bielefeld.de\\
X. Zhang: XichengZhang@gmail.com
 }

\begin{abstract}

In this paper, we prove the existence of a unique strong solution to a
stochastic tamed 3D Navier-Stokes equation in the whole space as well as in the periodic
boundary case.
Then, we also study the Feller property of solutions,
 and prove the existence of invariant measures
for the corresponding Feller semigroup in the case of periodic conditions.
Moreover, in the case of periodic boundary and degenerated additive noise, using the notion
of asymptotic strong Feller property proposed by Hairer and Mattingly
\cite{Ha-Ma}, we prove the uniqueness of invariant measures
for the corresponding transition semigroup.
\end{abstract}

\thanks{{\it AMS Classification(2000)}: 60H15, 37A25}

\maketitle \rm
\tableofcontents
\section{Introduction}

The classical 3D Navier-Stokes equations (NSE) describe the time evolution of an incompressible
fluid and are given by
$$
\p_t\u(t)=\nu\Delta \u(t)-(\u(t)\cdot\nabla)\u(t)+\nabla p(t)+\f(t)
$$
and
$$
\div \u(t)=0,
$$
where $\u(t,x)=(u^1(t,x),u^2(t,x),u^3(t,x))$ represents the velocity field,
$\nu$ is the viscosity constant,
$p(t,x)$ denotes the pressure, and $\f$ is an external force field acting on the fluid.
In \cite{Le}, Leray initially  constructed a weak solution for
the Cauchy problem of NSE in the whole space, since then,
it is still not known whether there exists a smooth  solution existing for all times.
In \cite{Ro-Zh}, we analyzed the following tamed scheme for the classical 3D NSE:
$$
\p_t\u(t)=\nu\Delta \u(t)-(\u(t)\cdot\nabla)\u(t)+\nabla p(t)-g_N(|\u(t)|^2)\u(t)+\f(t),
$$
where the taming function $g_N:\mR_+\to\mR_+$ is smooth and satisfies for some $N\in\mN$,
\begin{align}
\label{Con}\left\{
\begin{array}{ll}
g_N(r)=0, &\mbox{ if } r\leq N, \\
g_N(r)=(r-N)/\nu,& \mbox{ if } r\geq N+1,\\
0\leq g_N'(r)\leq 2/(\nu\wedge 1),& ~~r\geq 0.
\end{array}
\right.
\end{align}
Therein, we proved the existence of smooth solutions to this tamed equation when $\f$ and the initial velocity
are smooth. The main feature of this tamed equation is that if there is a bounded smooth solution
to the classical 3D NSE, then this smooth solution must satisfy our tamed equation
for some $N$ large enough. Moreover, we can let $N\to\infty$ to obtain
the existence of suitable weak solutions  (cf. \cite{Ro-Zh}).
In this sense, the above tamed scheme  can be considered as a regularized equation
for the classical equation.

Following the above tamed scheme, in the present paper we shall study the stochastic tamed 3D NSE.
Let us now describe our model equation.
Let $\mD:=\mR^3$ or $\mT^3$ (in the periodic case), where $\mT=[0,1)$ is the unit circle.
Note that any function from $\mT^3$ to $\mR^3$ can be identified with a periodic function
from $\mR^3$ to $\mR^3$.
We consider the following stochastic tamed 3D Navier-Stokes equation with $\nu=1$
in $\mD$:
\begin{align}
\dif \u(t)&=\Big[\Delta \u(t)-(\u(t)\cdot\nabla)\u(t)+\nabla p(t)-g_N(|\u(t)|^2)\u(t)
+\f(t,\u(t))\Big]\dif t\no\\
&\quad+\sum_{k=1}^\infty\Big[(\sigma_k(t)\cdot\nabla)\u(t)+\nabla \tilde p_k(t)
+\h_k(t,\u(t))\Big]\dif W^k_t\label{Ns1}
\end{align}
subject to the incompressibility condition
\begin{align}
\div\u(t)=0
\end{align}
and the initial condition
\begin{align}
\u(0)=\u_0,
\end{align}
where $p(t,x)$ and $\tilde p_k(t,x)$ are unknown scalar functions,
 $N>0$ and the taming function $g_N:\mR^+\to\mR^+$ as above satisfies (\ref{Con}),
and $\{W^k_t; t\geq 0, k=1,2,\cdots\}$ is a sequence of independent one dimensional
standard Brownian motions  on some complete filtration probability space
$(\Omega,\cF, P; (\cF_t)_{t\geq 0})$.
The stochastic integral is understood as It\^o's integral.
The entries of the coefficients are given as follows:
\begin{align*}
\mR_+\times\mD\times\mR^3\ni(t,x,\u)&\to\f(t,x,\u)\in\mR^3,\\
\mR_+\times\mD\ni(t,x)&\to\sigma(t,x)\in\mR^3\times l^2,\\
\mR_+\times\mD\times\mR^3\ni(t,x,\u)&\to\h(t,x,\u)\in\mR^3\times l^2,
\end{align*}
where $l^2$ denotes the Hilbert space consisting of all sequences of square summable
real numbers with standard norm $\|\cdot\|_{l^2}$.
In the following, $\f,\sigma$ and $\h$ are always assumed to be measurable
with respect to all their variables.

The study of stochastic Navier-Stokes equations (SNSE) began with
the work of Bensoussan -Temam in \cite{Be-Te}.
Using Galerkin's approximation and compactness method,
Flandoli  and Gatarek in \cite{Fl-Ga} proved the existence of martingale solutions
and stationary solutions for any dimensional
stochastic Navier-Stokes equations in a bounded domain. In particular, when the transition semigroup
is well defined, the stationary martingale solutions will yield the existence of invariant measures.
We remark that their results cannot be used in the case of whole space because of the
absence of  compact Sobolev embeddings.
Recently, Mikulevicius and Rozovskii in \cite{Mi-Ro} proved the existence of martingale solutions to
 SNSE in $\mR^d$ ($d\geq 2$) under less assumptions on the coefficients (without the extra term $g_N$).
To avoid the use of compact Sobolev embeddings, they used the approach of mollifying and cutting off the coefficients.
In the case of two dimension, they also obtained the existence and
pathwise uniqueness of $L^2$-continuous adapted solutions.

On the other hand, the ergodicity of invariant measures for 2D stochastic
Navier-Stokes equations has been studied extensively
 (cf. \cite{Fl-Ma,Ma,EMS,Ha-Ma} and reference therein). Especially, Hairer and Mattingly  \cite{Ha-Ma}
recently developed two important tools: the asymptotic strong Feller property and an approximative
integration by parts formula in the Malliavin calculus,
and then used them to derive an optimal ergodicity
result for 2D SNSE in the sense that the random forces only has two modes.
As pointed out in \cite{Ha-Ma}, the asymptotic strong Feller property is much weaker than
the usual strong Feller property since many degenerated equations have the former property
rather than the later one.

Up to now, to the best of  our knowledge, most of the well known results about the stochastic
Navier-Stokes equations such as the existence of invariant measures and
the ergodicity under different conditions on the noise are for 2D SNSE.
As for the three dimensional case,
there are only a few results  (cf. \cite{DaDe, DeOd, Ba, Od, Fl-Ro, Rom}),
of course, because of the lack of uniqueness.
Recently, in \cite{DaDe, DeOd, Od},
Da-Prato, Debussche and Odasso proved the existence and ergodicity of Markov solutions
 for 3D SNSE without the taming term $g_N$, which are obtained as limits of Galerkin's approximations.
Similar results were obtained by Flandoli and Romito
in \cite{Fl-Ro, Rom} for all Markov solutions. Moreover, using stochastic cascades,
Bakhtin \cite{Ba} explicitly constructed a stationary solution of 3D Navier-Stokes system
and proved a uniqueness theorem.

In the present paper, we shall prove the existence of a unique strong solution to our stochastic
tamed 3D Navier-Stokes equation (\ref{Ns1}) under some assumptions on $\f$, $\sigma$ and $\h$.
Here, the word ``strong'' means ``strong'' both in the sense of
the theory of stochastic differential equations
and the theory of partial differential equations. Let $\mH^m$ denote
the Sobolev space of divergence free vector fields (see (\ref{Div}) below).
Instead of working on the evolution
triple $\mH^1\subset\mH^0\subset\mH^{-1}$, we shall work on the evolution triple
$\mH^2\subset\mH^1\subset\mH^0$.
This will enable us to obtain the ``strong'' solution in the sense of
partial differential equations. For the ``strong'' solution in
the sense of stochastic differential equations, we shall use the famous Yamada-Watanabe theorem:
the existence of martingale solutions plus pathwise uniqueness implies the existence of a
unique strong solution.
Different from the method in \cite{Mi-Ro}, we still
use the classical Galerkin approximation to prove our existence of strong solutions.
To overcome the absence of compact Sobolev embeddings, we shall use localization method to prove
tightness. We think that it is of interest in itself and can be used in other cases.
Moreover, as in the deterministic case, we can take  limits $N\to\infty$
to prove the existence of weak solutions for the true stochastic Navier-Stokes equations
(without taming term). This will be done in a further investigation.

After obtaining the existence of a unique strong solution to Eq. (\ref{Ns1}),
we turn to the study of uniqueness of invariant measures in the case of
periodic boundary conditions and degenerated additive noise. As a first step, we need to
prove the Feller property and the existence of an invariant measure.
Then, using the asymptotic strong Feller property and approximative
integration by parts formula in \cite{Ha-Ma}, we can prove the uniqueness of
invariant measures. As said above, since we shall work in the first order Sobolev space $\mH^1$,
all of our discussions will take place in $\mH^1$. This requires
some delicate analysis and calculations, and the special form of $g_N$
plays an important role throughout this paper.
It should be emphasized that
the optimal results in \cite{Ha-Ma} seem to depend strongly on the structure of 2D
Navier-Stokes equations, we can not develop a similar non-adapted analysis along their lines
to obtain some optimal result for our tamed 3D SNSE.

This paper is organized as follows: in Section 2, we give some preliminaries, that include some
necessary estimates and a tightness result for later use.
In Section 3, we shall prove the existence and uniquess result by Galerkin's approximation.
In Section 4, we study the Feller property of the solutions to Eq. (\ref{Ns1}) and
the existence of invariant measures for the Feller semigroup
in the case of periodic boundary conditions.
In Section 5, we study the ergodicity of invariant measures. In the Appendix, for the reader's
convenience, the martingale characterization for weak solutions is proved,
two necessary basic estimates are given,  and the derivative flow equation
is proved.

\section{Preliminaries}
\subsection{Notations and Assumptions}
Let $C^\infty_0(\mD;\mR^3)$
denote the set of all smooth functions from $\mD$ to $\mR^3$ with
compact supports. When $\mD=\mT^3$, a function $f\in C^\infty_0(\mD;\mR^3)$
means that it is a smooth periodic function from $\mR^3$ to $\mR^3$.
For $p\geq 1$, let $L^p(\mD;\mR^3)$ be the vector valued $L^p$-space in which the norm
is denoted by $\|\cdot\|_{L^p}$.
For $m\in\mN_0:=\mN\cup\{0\}$,
let $H^m$ be the usual Sobolev space on $\mD$ with values in $\mR^3$,
i.e., the closure of
$C^\infty_0(\mD;\mR^3)$ with respect to the norm:
\begin{align*}
\|\u\|_{H^m}=\left(\int_{\mD}|(I-\Delta)^{m/2}\u|^2\dif x\right)^{1/2}.
\end{align*}
Here as usual, $(I-\Delta)^{m/2}$ is defined by Fourier transformation.
For two separable Hilbert spaces $\mK$ and $\mH$, $L_2(\mK;\mH)$ will denote the space of
all Hilbert-Schmidt operators from $\mK$ to $\mH$ with norm  $\|\cdot\|_{L_2(\mK;\mH)}$.

The following Gagliardo-Nirenberg interpolation inequality will be used frequently.
It plays an essential role in the study of Navier-Stokes equations (cf. \cite{Ta}).
Let $q\in[1,\infty]$ and $m\in\mN$. If
$$
\frac{1}{q}= \frac{1}{2}-\frac{m\a}{3},\quad 0\leq\a\leq 1,
$$
then for any $\u\in H^{m}$
\begin{align}
\|\u\|_{L^q}\leq C_{m,q}\|\u\|^\a_{H^m}\|\u\|_{L^2}^{1-\a}.\label{Sob}
\end{align}

Set for $m\in\mN_0$
\begin{align}
\mH^m:=\{\u\in H^{m}: \div(\u)=0\}.\label{Div}
\end{align}
Then $(\mH^m,\|\cdot\|_{H^m})$
is a  separable Hilbert space. We shall denote the norm
$\|\cdot\|_{H^m}$ in $\mH^m$ by $\|\cdot\|_{\mH^m}$.
We remark that $\mH^0$ is a closed linear subspace of the Hilbert space $L^2(\mD;\mR^3)=H^0$.

Let $\sP$ be the orthogonal projection  from $L^2(\mD;\mR^3)$ to $\mH^0$
 (cf.  \cite{Fa-Jo-Ri, Lions}).
It is well known that $\sP$ commutes with the derivative operators, and that
$\sP$ can be restricted to a bounded linear operator from $H^{m}$ to $\mH^m$.
For any $\u\in\mH^0$ and $\v\in L^2(\mD;\mR^3)$, we have
\begin{align*}
\<\u,\v\>_{\mH^0}:=\<\u,\sP\v\>_{\mH^0}=\<\u,\v\>_{L^2}.
\end{align*}

Let $\cV$ be defined by
\begin{align*}
\cV:=\{\u: \u\in C^\infty_0(\mD;\mR^3), \div(\u)=0\}.
\end{align*}
We have the following density result (cf. \cite{Ro-Zh}).
\bl\label{Le3}
$\cV$ is dense in $\mH^m$ for any $m\in\mN_0$.
\el

We now introduce the following assumptions on the coefficients $\f,\sigma$ and $\h$:
\begin{enumerate}[{\bf (H1)}]
\item For any $T>0$, there exist a constant $C_{T,\f}>0$ and a
 function $H_{\f}(t,x)\in L^1([0,T]\times\mD)$ such that for any $t\in[0,T], x\in\mD,\u\in\mR^3$
and $j=1,2,3$
\begin{align*}
|\p_{x^j}\f(t,x,\u)|^2+|\f(t,x,\u)|^2&\leq C_{T,\f}\cdot|\u|^2+H_{\f}(t,x),\\
|\p_{u^j}\f(t,x,\u)|&\leq C_{T,\f}.
\end{align*}
\item For any $T>0$, there exists a constant $C_{\sigma,T}>0$ such that
$$
\sup_{t\in[0,T],x\in\mD}\|\p_{x_j}\sigma(t,x)\|_{l^2}\leq C_{\sigma,T},\ \ j=1,2,3
$$
and
\begin{align}
\sup_{t\in\mR_+,x\in\mD}\|\sigma(t,x)\|_{l^2}^2\leq 1/4.\label{Sigma}
\end{align}
\item For any $T>0$, there exist a constant $C_{T,\h}>0$ and a
 function $H_{\h}(t,x)\in L^1([0,T]\times\mD)$ such that for any $t\in[0,T], x\in\mD,\u\in\mR^3$
and $j=1,2,3$
\begin{align*}
\|\p_{x^j}\h(t,x,\u)\|_{l^2}^2+\|\h(t,x,\u)\|_{l^2}^2&\leq C_{T,\h}\cdot|\u|^2+H_{\h}(t,x),\\
\|\p_{u^j}\h(t,x,\u)\|_{l^2}&\leq C_{T,\h}.
\end{align*}
\end{enumerate}

\br
The factor $\frac{1}{4}$ in (\ref{Sigma}) is related to the viscosity constant $\nu$
assumed to be $1$. That is to say, the first order term appearing in
diffusion coefficients will
be absorbed by the Laplace term. Here, the factor $\frac{1}{4}$ is not optimal (see \cite{Mi-Ro}).
\er

For any $\u\in\mH^2$, define
\begin{align}
A(\u):=\sP\Delta\u-\sP((\u\cdot\nabla)\u)-\sP(g_N(|\u|^2)\u),\label{AA}
\end{align}
and for any $\v\in\mH^2$, we write
$$
\lb A(\u),\v\rb:= \<A(\u),(I-\Delta)\v\>_{\mH^0}= A_1(\u,\v)+A_2(\u,\v)+A_3(\u,\v),
$$
where
\begin{align*}
A_1(\u,\v)&:= \<\Delta \u,(I-\Delta) \v\>_{\mH^0},\\
A_2(\u,\v)&:= -\<(\u\cdot\nabla)\u, (I-\Delta)\v\>_{\mH^0},\\
A_3(\u,\v)&:= -\<g_N(|\u|^2)\u, (I-\Delta)\v\>_{\mH^0}.
\end{align*}

Below, for the sake of simplicity,   the variable
``$x$'' in the coefficients will be dropped.
Define for $k\in\mN$
\begin{align}
B_k(t,\u):=\sP((\sigma_k(t)\cdot\nabla)\u)+\sP\h_k(t,\u).\label{BB}
\end{align}

Letting the operator $\sP$ act on both sides of equation (\ref{Ns1}), we can and shall
consider the following equivalent abstract stochastic evolution equation in the sequel:
\begin{align}
\label{NS}\left\{
\begin{array}{lcl}
\dif \u(t)&= &\Big[A(\u(t))+\sP\f(t,\u(t))\Big]\dif t+\sum_{k=1}^\infty B_k(t,\u(t))\dif W^k_t,\\
\u(0)&=& \u_0\in\mH^1.
\end{array}
\right.
\end{align}

\subsection{Estimates on $A$ and $B$}

We now prepare several important estimates for later use.
In the sequel, we shall use the following convention:  The letter $C$ with subscripts will
denote a constant depending on its subscripts and the coefficients. The letter $C$ without subscripts
will denote an absolute constant, i.e., its value
does not depend on any data. All the constants may have
different values in different places.

\bl\label{Le1}
For any $\u\in\mH^2$, we have
\begin{align}
\|A(\u)\|_{\mH^0}&\leq C(1+\|\u\|^4_{\mH^0}+\|\u\|^2_{\mH^2}),\label{Op8}\\
\<A(\u),\u\>_{\mH^0}&= -\|\nabla\u\|^2_{\mH^0}-\|\sqrt{g_N(|\u|^2)}\cdot|\u|\|^2_{L^2}\label{Es444}\\
&\leq -\|\nabla\u\|^2_{\mH^0}-\|\u\|^4_{L^4}+C\cdot N\|\u\|^2_{\mH^0},\label{Es44}\\
\lb A(\u),\u\rb&\leq -\frac{1}{2}\| \u\|^2_{\mH^2}-\frac{1}{2}\||\u|\cdot|\nabla\u|\|^2_{L^2}+
C\cdot N\|\nabla\u\|^2_{\mH^0}+\|\u\|^2_{\mH^0}.\label{Es4}
\end{align}
\el
\begin{proof}
Estimate (\ref{Op8}) is direct from (\ref{AA}) and
the Sobolev inequality (\ref{Sob}). Estimate (\ref{Es444}) follows from
$$
\<(\u\cdot\nabla)\u,\u\>_{\mH^0}=0.
$$
For inequality (\ref{Es4}), we have
\begin{align*}
A_1(\u,\u)= -\|(I-\Delta)\u\|^2_{\mH^0}+\<\u, (I-\Delta)\u\>_{\mH^0}
= -\|\u\|^2_{\mH^2}+\|\nabla\u\|^2_{\mH^0}+\|\u\|^2_{\mH^0},
\end{align*}
and by Young's inequality,
\begin{align*}
A_2(\u,\u)\leq \frac{1}{2}\|(I-\Delta)\u\|^2_{\mH^0}
+\frac{1}{2}\|(\u\cdot\nabla)\u\|^2_{\mH^0}
\leq \frac{1}{2}\|\u\|^2_{\mH^2}+\frac{1}{2}\||\u|\cdot|\nabla\u|\|^2_{\mH^0},
\end{align*}
where
$$
|\u|^2=\sum_{k=1}^3|u^k|^2,\quad |\nabla\u|^2=\sum_{k,i=1}^3|\p_iu^k|^2.
$$
Recalling $\nu=1$, from (\ref{Con}), we  also have
\begin{align*}
A_3(\u,\u)&= -\<\nabla (g_N(|\u|^2)\u), \nabla\u\>_{\mH^0}-\<g_N(|\u|^2)\u, \u\>_{\mH^0}\\
&= -\sum_{k,i=1}^3\int_{\mD} \p_i u^k\cdot \p_i (g_N(|\u|^2)u^k)\dif x
-\int_{\mD}|\u|^2\cdot g_N(|\u|^2)\dif x\\
&\leq -\sum_{k,i=1}^3\int_{\mD} \p_i u^k\cdot \left(g_N(|\u|^2)\cdot \p_i u^k
-g'_N(|\u|^2)\p_i |\u|^2\cdot u^k\right)\dif x\\\
&= -\int_{\mD}|\nabla \u|^2\cdot g_N(|\u|^2)\dif x
-\frac{1}{2}\int_{\mD}g'_N(|\u|^2)|\nabla|\u|^2|^2\dif x\\
&\leq -\int_{\mD}|\nabla \u|^2\cdot |\u|^2\dif x+C\cdot N\|\nabla\u\|^2_{\mH^0}.
\end{align*}
Combining the above calculations yields (\ref{Es4}).
\end{proof}

\bl\label{Le2}
Let $\v\in\cV$, and let the support of $\v$ be contained in $\cO:=\{x\in\mD, |x|\leq m\}$ for some $m\in\mN$.
Let $T>0$. For any $\u,\u'\in\mH^2$ and $t\in[0,T]$, we have
\begin{align}
|\lb A(\u),\v\rb|&\leq C_{\v}\cdot\big(1+\|\u\|^3_{L^3(\cO)}\big),\label{Op1}\\
\|\<B_\cdot(t,\u),\v\>_{\mH^1}\|^2_{l^2}&\leq C_{\v,T}
\cdot\big(1+\|H_{\h}(t)\|_{L^1(\mD)}+\|\u\|^2_{L^2(\cO)}\big)\label{Op2}
\end{align}
and
\begin{align}
|\lb A(\u)-A(\u'),\v\rb|\leq C_{\v}\cdot\|\u-\u'\|_{L^2(\cO)}\cdot(1+\|\u\|^2_{\mH^1}
+\|\u'\|^2_{\mH^1}).\label{Op3}
\end{align}
\el
\begin{proof}
For estimate (\ref{Op1}), we have
\begin{align*}
A_1(\u,\v)=\<\u, (I-\Delta)\Delta \v\>_{\mH^0}&\leq C\|\u\|_{L^2(\cO)}\cdot\|\v\|_{\mH^4},\\
A_2(\u,\v)=\<\u^*\cdot\u,\nabla(I-\Delta)\v\>_{\mH^0}
&\leq C\|\u\|^2_{L^2(\cO)}\cdot\sup_{x\in\mD}|\nabla(I-\Delta) \v(x)|,
\end{align*}
where $\u^*$ denotes the transposition of the row vector $\u$,
and
$$
A_3(\u,\v)\leq \|\u\|^3_{L^3(\cO)}\cdot\sup_{x\in\mD}|(I-\Delta) \v(x)|.
$$
Combining them gives (\ref{Op1}).

For estimate (\ref{Op2}),  by {\bf (H2)} and {\bf (H3)}, we have
\begin{align*}
\|\<B_\cdot(t,\u),\v\>_{\mH^1}\|^2_{l^2}&\leq C\sup_{t\in\mR_+,x\in\mD}\|\sigma(t,x)\|^2_{l^2}
\cdot\sup_{x\in\mD}|\nabla(I-\Delta)\v(x)|^2\cdot\|\u\|^2_{L^2(\cO)}\\
&\quad +C\sup_{x\in\mD}\|\nabla_x\sigma(t,x)\|^2_{l^2}
\cdot\sup_{x\in\mD}|(I-\Delta)\v(x)|^2\cdot\|\u\|^2_{L^2(\cO)}\\
&\quad +C\sup_{x\in\mD}|(I-\Delta)\v(x)|^2\cdot\big(C_{T,\h}\cdot\|\u\|^2_{L^2(\cO)}
+\|H_{\h}(t)\|_{L^1(\mD)}\big)\\
&\leq C_{\v,T}\cdot\big(1+\|H_{\h}(t)\|_{L^1(\mD)}+\|\u\|^2_{L^2(\cO)}\big).
\end{align*}

We now look at (\ref{Op3}).
For $A_1$, we clearly have
\begin{align*}
|A_1(\u,\v)-A_1(\u',\v)|
=|\<(\u-\u')\cdot 1_\cO, (I-\Delta)\Delta \v\>_{\mH^0}|\leq C_{\v}\cdot\|\u-\u'\|_{L^2(\cO)}.
\end{align*}
For $A_2$, we have
\begin{align*}
|A_2(\u,\v)-A_2(\u',\v)|
&= |\<\u^*\cdot\u-\u'^*\cdot\u', \nabla(I-\Delta) \v\>_{\mH^0}|\\
&\leq C_{\v}\cdot\|\u-\u'\|_{L^2(\cO)}\cdot(\|\u\|_{\mH^0}+\|\u'\|_{\mH^0}).
\end{align*}
For $A_3$, by Sobolev  inequality (\ref{Sob}) we similarly have
\begin{align*}
|A_3(\u,\v)-A_3(\u',\v)|
&\leq C_{\v}\cdot\|\u-\u'\|_{L^2(\cO)}\cdot(\|\u\|_{L^4}^2+\|\u'\|_{L^4}^2)\\
&\leq C_{\v}\cdot\|\u-\u'\|_{L^2(\cO)}\cdot(\|\u\|_{\mH^1}^2+\|\u'\|_{\mH^1}^2).
\end{align*}
\end{proof}

\bl\label{Le6}
For  any $T>0$ and $\u\in\mH^2$,
\begin{align}
\|B(t,\u)\|^2_{L_2(l^2;\mH^0)}&\leq \frac{1}{2}\|\u\|^2_{\mH^1}+C_T\|\u\|^2_{\mH^0}
+\|H_{\h}(t)\|_{L^1(\mD)},\label{LP1}\\
\|B(t,\u)\|^2_{L_2(l^2;\mH^1)}&\leq \frac{1}{2}\|\u\|^2_{\mH^2}+C_T\|\u\|^2_{\mH^1}
+C\|H_{\h}(t)\|_{L^1(\mD)}.\label{LP2}
\end{align}
\el
\begin{proof}
First of all, by {\bf (H2)} and {\bf (H3)}, we have
\begin{align*}
& \|B(t,\u)\|^2_{L_2(l^2;\mH^0)}=\sum_{k=1}^\infty\int_{\mD}|B_k(t,x,\u(x))|^2\dif x\leq\\
&\qquad\leq 2\int_{\mD}\|\sigma(t,x)\|_{l^2}^2\cdot |\nabla\u(x)|^2\dif x
+2\int_{\mD}\big(C_{T,\h}|\u(x)|^2+H_{\h}(t,x)\big)\dif x\\
&\qquad\leq \frac{1}{2}\|\u\|^2_{\mH^1}+C_T\|\u\|^2_{\mH^0}
+\|H_{\h}(t)\|_{L^1(\mD)}.
\end{align*}
Secondly, noting that
$$
\|B(t,\u)\|^2_{L_2(l^2;\mH^1)}=\|B(t,\u)\|^2_{L_2(l^2;\mH^0)}+\|\nabla B(t,\u)\|^2_{L_2(l^2;\mH^0)}
$$
and
\begin{align*}
\p_{x^j}B_k(t,\u)&= \sP\p_{x^j}((\sigma_k(t)\cdot\nabla)\u)+\sP\p_{x^j}\h_k(t,\u)\\
&= \sP\Big((\p_{x^j}\sigma_k(t)\cdot\nabla)\u+(\sigma_k(t)\cdot\nabla)\p_{x^j}\u\Big)\\
&\quad +\sP\Big((\p_{x^j}\h_k)(t,\u)+\sum_{i=1}^3\p_{u^i}\h_k(t,\u)\cdot\p_{x^j}u^i\Big),
\end{align*}
by {\bf (H2)} and {\bf (H3)}, we have
\begin{align*}
\|B(t,\u)\|^2_{L_2(l^2;\mH^1)}\leq \frac{1}{2}\|\u\|^2_{\mH^2}+
C_T\|\u\|^2_{\mH^1}+C\|H_{\h}(t)\|_{L^1(\mD)}.
\end{align*}
\end{proof}

\subsection{Tightness Criterion}

In the following, we only give a tightness criterion in the case of $\mD=\mR^3$.
When $\mD=\mT^3$, since $\mH^1$ is compactly embedded in $\mH^0$,
the corresponding result is simple and well known.

By $\mH^0_{loc}$ we denote the space of all locally $L^2$-integrable and
divergence free vector fields  endowed
with the Fr\'echet metric:  for $\u,\v\in\mH^0_{loc}$
$$
\rho(\u,\v):=\sum^\infty_{m=1}2^{-m}\left(\left[\int_{|x|\leq m}
|\u(x)-\v(x)|^2\dif x\right]^{1/2}\wedge 1\right).
$$
Thus, $(\mH^0_{loc},\rho)$ is a Polish space and $\mH^0\subset\mH^0_{loc}$.

Let $\mX:=C(\mR_+;\mH^0_{loc})$ denote the space of all continuous functions from $\mR_+$
to $(\mH^0_{loc},\rho)$ equipped with the metric
$$
\rho_\mX(\u,\v):=\sum^\infty_{m=1}2^{-m}\left(\sup_{t\in[0,m]}\rho(\u(t),\v(t))\wedge 1\right).
$$

In the following, we shall fix a complete orthonormal basis $\sE:=\{\e_k,k\in\mN\}\subset\cV$ of $\mH^1$
such that $\mbox{span}\{\sE\}$ is a dense subset of $\mH^3$ and, in the case of periodic boundary
conditions,
we  also require that $\sE$ is an orthogonal basis of $\mH^0$.
Moreover, for $\u\in\mH^0$ and $\v\in\mH^2$,  the inner product
$\<\u,\v\>_{\mH^1}$ is taken in the generalized sense,
i.e.,
$$
\<\u,\v\>_{\mH^1}=\<\u,(I-\Delta)\v\>_{\mH^0}.
$$

We need the following relative  compactness result, which is essentially due to
Ladyzhenskaya \cite[Theorem 13]{La}.
\bl\label{Le4}
Let $K\subset \mX$. If for every $T>0$,
\begin{enumerate}[($1^o$)]
\item $\sup_{\u\in K}\sup_{s\in[0,T]}\|\u(s)\|_{\mH^1}<+\infty$,
\item $\lim_{\delta\rightarrow 0}\sup_{\u\in K}\sup_{t,s\in[0,T],|t-s|<\delta}
|\<\u(t)-\u(s),\e\>_{\mH^1}|=0$ for any $\e\in\sE$,
\end{enumerate}
then $K$ is relatively compact in $\mX$.
\el
\begin{proof}
We only need to prove that $K$ is relatively compact in $C([0,T];\mH^0_{loc})$ for every $T>0$.
Let $\{\u_n,n\in\mN\}\subset K$ be any sequence of $K$. Define for $\e\in\sE$,
\begin{align*}
G^{\e}_n(t):=\<\u_n(t),\e\>_{\mH^1}.
\end{align*}
Then, by ($1^o$) and ($2^o$),
the sequence $\{t\mapsto G^{\e}_n(t), n\in\mN\}$ is uniformly bounded and equi-continuous on $[0,T]$.
Hence, by Ascoli-Arzel\`a's lemma, there exist a subsequence
$n_l$ (depending on $\e$) and a continuous function $G^{\e}(t)$ such that
$G^{\e}_{n_l}(t)$ uniformly converges to
$G^{\e}(t)$ on $[0,T]$. Since $\sE$ is countable,
by a diagonalization method, we may further find a common subsequence
(still denoted by $n$) such that for any $\e\in\sE$,
$$
\lim_{n\rightarrow\infty}\sup_{t\in[0,T]}|G^{\e}_n(t)-G^{\e}(t)|=0.
$$
Thus, by the weak compactness of closed balls in $\mH^1$,
there is a $\u\in L^\infty(0,T;\mH^1)$ such that for any $\e\in\sE$,
$$
\lim_{n\rightarrow\infty}\sup_{t\in[0,T]}|\<\u_n(t)-\u(t),\e\>_{\mH^1}|=0.
$$
By a simple approximation we further have for any $\v\in\mH^1$,
$$
\lim_{n\rightarrow\infty}\sup_{t\in[0,T]}|\<\u_n(t)-\u(t),\v\>_{\mH^1}|=0.
$$
Note that $(I-\Delta)^{-1}\v\in\mH^2$ for any $\v\in\mH^0$. Hence, we also have for any $\v\in\mH^0$,
$$
\lim_{n\rightarrow\infty}\sup_{t\in[0,T]}|\<\u_n(t)-\u(t),\v\>_{\mH^0}|=0.
$$
Hence, by Helmholtz-Weyl's decomposition (cf. \cite{Te, Ga}),
\begin{align}
\lim_{n\rightarrow\infty}\sup_{t\in[0,T]}|\<\u_n(t)-\u(t),\v\>|=0\label{WEAK}
\end{align}
for any $\v\in L^2(\mR^3;\mR^3)$.

We now show that
\begin{align*}
\lim_{n\rightarrow\infty}\sup_{t\in[0,T]}\rho(\u_n(t),\u(t))=0.
\end{align*}
It suffices to prove that for any $m\in\mN$,
\begin{align*}
\lim_{n\rightarrow\infty}\sup_{t\in[0,T]}\int_{|x|\leq m}|\u_n(t,x)-\u(t,x)|^2\dif x=0,
\end{align*}
which follows from ($1^o$), (\ref{WEAK}) and the following Friedrichs inequality
 (cf. \cite[p.176]{La}):
Let $\cO\subset\mR^3$ be any bounded  domain. For any
$\epsilon>0$, there exist $N_\epsilon\in\mN$ and
functions $h_i\in L^2(\cO), i=1,\cdots, N_\epsilon$ such that for any $w\in W^{1,2}_0(\cO)$,
\begin{align*}
\int_\cO|w(x)|^2\dif x\leq
\sum^{N_\epsilon}_{i=1}\left(\int_\cO w(x) h_i(x)\dif x\right)^2
+\epsilon\int_\cO|\nabla w(x)|^2\dif x.
\end{align*}
\end{proof}

\bl\label{Le5}
Let $\mu_n$ be a family of probability measures on $(\mX, \cB(\mX))$. Assume that
\begin{enumerate}[($1^o$)]

\item For each $\e\in\sE$ and any $\epsilon, T>0$,
$$
\lim_{\delta\downarrow 0}\sup_n\mu_n
\Big\{\u\in\mX: \sup_{ s,t\in[0,T],|s-t|\leq \delta}|\<\u(t)-\u(s),\e\>_{\mH^1}|>\epsilon\Big\}=0.
$$

\item For any $T>0$
$$
\lim_{R\rightarrow\infty}\sup_n \mu_n\Big\{\u\in\mX: \sup_{s\in[0,T]}\|\u(s)\|_{\mH^1}>R\Big\}=0.
$$
\end{enumerate}
Then $\{\mu_n,n\in\mN\}$ is tight on $(\mX, \cB(\mX))$.
\el
\begin{proof}
Fix $\eta>0$.  For any $l\in\mN$, by ($2^o$) one can choose $R_l$ sufficiently large such that
\begin{align}
\sup_n \mu_n\Big\{\u\in\mX: \sup_{s\in[0,l]}\|\u(s)\|_{\mH^1}> R_l\Big\}\leq \frac{\eta}{2^l}.\label{Es3}
\end{align}
For $k,l\in\mN$ and $\e_i\in\sE$, by ($1^o$) one may choose $\delta_{k,i,l}>0$ small enough such that
\begin{align}
\sup_n\mu_n
\Big\{\u\in\mX: \sup_{s,t\in[0,l],|s-t|\leq \delta_{k,i,l}}|\<\u(t)-\u(s),\e_i\>_{\mH^1}|
>\frac{1}{k}\Big\}\leq\frac{\eta}{2^{k+i+l}}.\label{Es2}
\end{align}

Now let us define
\begin{align*}
K_1&:= \bigcap_{k,l\in\mN,\e_i\in\sE}\Big\{\u\in\mX: \sup_{s,t\in[0,l],|s-t|\leq \delta_{k,i,l}}|\<\u(t)-\u(s),
\e_i\>_{\mH^1}|\leq \frac{1}{k}\Big\}\\
K_2&:= \bigcap_{l\in\mN}\Big\{\u\in\mX: \sup_{s\in[0,l]}\|\u(s)\|_{\mH^1}\leq R_l\Big\}.
\end{align*}
By Lemma \ref{Le4}, $K_1\cap K_2$ is a relatively compact set in $\mX$. By (\ref{Es3}) and (\ref{Es2}), we also have
\begin{align*}
\sup_n\mu_n(K_1^c\cup K_2^c)\leq 2\eta.
\end{align*}
In view of the arbitrariness of $\eta$, $\{\mu_n,n\in\mN\}$ is tight on $(\mX, \cB(\mX))$.
\end{proof}

\section{Existence and Uniqueness of Strong Solutions}

\subsection{Weak and Strong Solutions}
For a metric space $\mU$, we use $\cP(\mU)$ to denote the total of all probability measures on $\mU$.
We first introduce the following notion of  weak solutions to Eq. (\ref{NS}).
\bd\label{def}
We say that Eq. (\ref{NS}) has a weak solution with initial law $\vartheta\in\cP(\mH^1)$ if there exist a
stochastic basis $(\Omega,\cF,P; (\cF_t)_{t\geq 0})$,  an $\mH^1$-valued $(\cF_t)$-adapted
process $\u$
and an infinite sequence of independent standard $(\cF_t)$-Brownian motions $\{W^k(t), t\geq 0, k\in\mN\}$
such that
\begin{enumerate}[(i)]
\item $\u(0)$ has law $\vartheta$ in $\mH^1$;
\item for  almost all $\om\in\Omega$ and every $T>0$,
$\u(\cdot,\om)\in C([0,T];\mH^1)\cap L^2([0,T];\mH^2)$;

\item it holds that in $\mH^0$
\begin{align*}
\u(t)=\u_0+\int^t_0\Big[A(\u(s))+\sP\f(s,\u(s))\Big]\dif s+\sum_{k=1}^\infty\int^t_0 B_k(s,\u(s))\dif W^k_s,
\end{align*}
for all $t\geq 0$, $P$-a.s..
\end{enumerate}

This solution is denoted by $(\Omega,\cF,P; (\cF_t)_{t\geq 0}; W; \u)$.
\ed

\br
Under {{\bf (H1)}}-{\bf (H3)}, by (\ref{Op8}) the above integrals are meaningful.
\er
\bd\label{def4}
(Pathwise Uniqueness) We say that the pathwise uniqueness holds  for Eq. (\ref{NS}) if
whenever we are given two  weak solutions of Eq. (\ref{NS})
defined on the same probability space together with
the same Brownian motion
\begin{align*}
&\quad (\Omega,\cF,P; (\cF_t)_{t\geq 0}; W; \u)\\
&\quad (\Omega,\cF,P; (\cF_t)_{t\geq 0}; W; \tilde \u),
\end{align*}
the condition $P\{\u(0)=\tilde \u(0)\}=1$ implies
$P\{\om: \u(t,\om)=\tilde \u(t,\om), \forall t\geq 0\}=1$.
\ed

We have the following martingale characterization for the weak solution  (cf. \cite{St}).
For the reader's convenience, a short proof is provided in the Appendix.
\bp\label{Th2}
Let $\sE$ be given in Subsection 2.3.
For $\vartheta\in\cP(\mH^1)$, the following two statements are equivalent:
\begin{enumerate}[(i)]
\item Eq. (\ref{NS}) has a weak solution with initial law $\vartheta$.

\item There exists a probability measure $P_\vartheta\in\cP(\mX)$ such that for $P_\vartheta$-almost all $\u\in\mX$
and any $T>0$,
\begin{align}
\u\in L^\infty([0,T];\mH^1)\cap L^2([0,T];\mH^2),\label{Es00}
\end{align}
and for any $h\in C_0^\infty(\mR)$, i.e., any smooth function with compact support, and any $\e\in\sE$,
\begin{align*}
M^h_{\e}(t,\u)&:= h(\<\u(t),\e\>_{\mH^1})-h(\<\u(0),\e\>_{\mH^1})\\
&\quad -\int^t_0h'(\<\u(s),\e\>_{\mH^1})\cdot\lb  A(\u(s)),\e\rb\dif s\\
&\quad -\int^t_0h'(\<\u(s),\e\>_{\mH^1})\cdot \<\f(s,\u(s)),\e\>_{\mH^1}\dif s\\
&\quad -\frac{1}{2}\int^t_0h''(\<\u(s),\e\>_{\mH^1})\cdot\|\<B(s,\u(s)),\e\>_{\mH^1}\|^2_{l^2}\dif s
\end{align*}
is a continuous local martingale under $P_\vartheta$ with respect to $\cB_t(\mX)$.
Here and below, $\cB_t(\mX)$ denotes the sub $\sigma$-algebra of $\mX$ up to time $t$.
\end{enumerate}
\ep

In order to introduce the notion of  strong solutions to Eq. (\ref{NS}),
we need a canonical realization of
an infinite sequence of independent standard Brownian motions on a Polish space.

Let $C(\mR_+;\mR)$ denote the space of all continuous functions
defined on $\mR_+$, which is equipped with the metric
$$
\tilde\rho(w,w')=\sum\limits_{k=1}^{\infty}2^{-k}\left(\sup_{t\in[0,k]}|w(t)-w'(t)|\wedge1\right).
$$
Define the product space $\mW:=\prod\limits_{j=1}^{\infty}C(\mR_+;\mR)$, which is endowed with
the metric:
$$
\rho_{\mW}(w,w^{\prime})=\sum\limits_{j=1}^{\infty}2^{-j}(\tilde\rho(w^j,w'^j)\wedge1),\quad w=(w^1,w^2,\cdots),
w'=(w'^1,w'^2,\cdots).
$$
Then $(\mW,\rho_\mW)$ is a Polish space.
Let $\cB_t(\mW)\subset\cB(\mW)$ be the $\sigma$-algebra up to time $t$.
We endow $(\mW,\cB(\mW))$ with the Wiener measure $\mP$ such that the coordinate process
$$
w(t):=(w^1(t),w^2(t),\cdots)
$$
is an infinite sequence of independent standard $\cB_t(\mW)$-Brownian motions on $(\mW,\cB(\mW),\mP)$.

Let $\mB:=C(\mR_+;\mH^1)$ denote the space of all continuous functions from $\mR_+$ to $\mH^1$,
which is endowed with the metric
$$
\rho_\mB(\u,\v):=\sum^\infty_{k=1}2^{-k}\left(\sup_{t\in[0,k]}\|\u(t)-\v(t)\|_{\mH^1}\wedge 1\right).
$$
In the following, $\cB_t(\mB)$ denotes the sub $\sigma$-algebra of $\mB$ up to time $t$.
For a measure space $(S,\cS,\lambda)$, $\overline{\cS}^\lambda$
will denote the completion of $\cS$ with respect to $\lambda$.

\bd\label{Def1}
Let $(\Omega,\cF,P; (\cF_t)_{t\geq 0}; W; \u)$ be a weak solution of Eq. (\ref{NS}) with initial
distribution $\vartheta\in\cP(\mH^1)$. If there exists a
$\overline{\cB(\mH^1)\times\cB(\mW)}^{\vartheta\times \mP}/\cB(\mB)$-measurable
functional $F_\vartheta:\mH^1\times\mW\mapsto \mB$, with the property that for every $t>0$,
\begin{align}
F_\vartheta\in\hat\cB_t/\cB_t(\mB); \quad \hat\cB_t:=
\overline{\cB(\mH^1)\times\cB_t(\mW)}^{\vartheta\times \mP}\label{P0}
\end{align}
and such that
\begin{align*}
\u(\cdot)=F_\vartheta(\u(0),W(\cdot)), \quad P-a.s.,
\end{align*}
we call $\u$ together with $W$ a strong solution.

We shall say that Eq. (\ref{NS}) has a unique strong solution associated with $\vartheta\in\cP(\mH^1)$
if there exists a functional $F_\vartheta:\mH^1\times\mW\mapsto \mB$ with the same properties
as above such that
\begin{enumerate}[(i)]
\item for any infinite sequence of independent
standard $(\cF_t)$-Brownian motions $\{W(t), t\geq 0\}$ on stochastic basis
$(\Omega,\cF,P; (\cF_t)_{t\geq 0})$, and any $\mH^1$-valued random
variable $\u_0\in\cF_0$ with distribution $\vartheta$,
$$
(\Omega,\cF,P; (\cF_t)_{t\geq 0}; W; F_\vartheta(\u_0,W(\cdot)))
\mbox{ is a weak solution of Eq. (\ref{NS})};
$$

\item for any weak solution $(\Omega,\cF,P; (\cF_t)_{t\geq 0}; W; \u)$ of Eq. (\ref{NS}) with initial law $\vartheta$,
\begin{align*}
\u(\cdot)=F_\vartheta(\u(0),W(\cdot)), \quad P-a.s..
\end{align*}
\end{enumerate}
\ed

The following Yamada-Watanabe theorem holds in this case  (cf. \cite{Ro-Zh1}).
\bt\label{Ya-Wa}
Existence of weak solutions plus pathwise uniqueness implies the existence
of a unique strong solution.
\et

\subsection{Pathwise Uniqueness}
We first prove the following pathwise uniqueness result.
\bt\label{Uni}
Under {{\bf (H1)}}-{\bf (H3)}, pathwise uniqueness holds for Eq. (\ref{NS}).
\et
\begin{proof}
Let $\u$ and $\tilde \u$ be two weak solutions of Eq. (\ref{NS})
defined on the same probability space together with the same Brownian motion, and starting from the
same initial value $\u_0$.
For any $T>0$ and $R>0$, define the stopping time
\begin{align*}
\tau_R:=\inf\{t\in[0,T]: \|\u(t)\|_{\mH^1}\vee\|\tilde \u(t)\|_{\mH^1}\geq R\}.
\end{align*}
By the definition of weak solutions, one knows that $\tau_R\uparrow \infty$ as $R\uparrow\infty$.

Set
$$
\w(t):=\u(t)-\tilde\u(t).
$$
Then by It\^o's formula, we have
\begin{align}
\|\w(t)\|^2_{\mH^0}&= 2\int^t_0\<A(\u(s))-A(\tilde\u(s)),\w(s)\>_{\mH^0}\dif s\no\\
&\quad +2\int^t_0\<\f(s,\u(s))-\f(s,\tilde\u(s)),\w(s)\>_{\mH^0}\dif s\no\\
&\quad +2\sum_{k=1}^\infty\int^t_0\<B_k(s,\u(s))-B_k(s,\tilde\u(s)),\w(s)\>_{\mH^0}\dif W^k_s\no\\
&\quad +\sum_{k=1}^\infty\int^t_0\|B_k(s,\u(s))-B_k(s,\tilde\u(s))\|^2_{\mH^0}\dif s\no\\
&=:I_1(t)+I_2(t)+I_3(t)+I_4(t).\label{Lp7}
\end{align}
By $|g_N(r)-g_N(r')|\leq |r-r'|$ and a simple calculation, it is easy to see that
\begin{align*}
I_1(t)&= -2\int^t_0\|\nabla\w(s)\|^2_{\mH^0}\dif s
+2\int^t_0\<\nabla\w(s),(\u^*(s)\cdot\u(s)-\tilde\u^*(s)\cdot\tilde\u(s))\>_{\mH^0}\dif s\\
&\quad -2\int^t_0\<g_N(|\u(s)|^2)\u(s)-g_N(|\tilde\u(s)|^2)\tilde\u(s),\w(s)\>_{\mH^0}\dif s\\
&\leq -\int^t_0\|\nabla\w(s)\|^2_{\mH^0}\dif s
+\int^t_0\|\u^*(s)\cdot\u(s)-\tilde\u^*(s)\cdot\tilde\u(s)\|^2_{\mH^0}\dif s\\
&\quad +8\int^t_0\||\w(s)|\cdot(|\u(s)|+|\tilde\u(s)|)\|^2_{\mH^0}\dif s.
\end{align*}
Noting that by Sobolev inequality (\ref{Sob}),
\begin{align}
\|\u^*(s)\cdot\u(s)-\tilde\u^*(s)\cdot\tilde\u(s)\|^2_{\mH^0}
&\leq \||\w(s)|(|\u(s)|+|\tilde\u(s)|)\|^2_{\mH^0}\label{Op00}\\
&\leq 2\|\w(s)\|^2_{L^4}(\|\u(s)\|^2_{L^4}+\|\tilde\u(s)\|^2_{L^4})\no\\
&\leq 2C_{1,4}^2\cdot\|\w(s)\|^{3/2}_{\mH^1}\|\w(s)\|^{1/2}_{\mH^0}
(\|\u(s)\|^2_{\mH^1}+\|\tilde\u(s)\|^2_{\mH^1}),\no
\end{align}
we have by Young's inequality,
\begin{align*}
I_1(t\wedge\tau_R)&\leq -\int^{t\wedge\tau_R}_0\|\nabla\w(s)\|^2_{\mH^0}\dif s
+C_R\int^{t\wedge\tau_R}_0\|\w(s)\|^{3/2}_{\mH^1}\|\w(s)\|^{1/2}_{\mH^0}\dif s\\
&\leq -\frac{1}{2}\int^{t\wedge\tau_R}_0\|\nabla\w(s)\|^2_{\mH^0}\dif s
+C_R\int^{t\wedge\tau_R}_0\|\w(s)\|^2_{\mH^0}\dif s.
\end{align*}
Moreover, it is clear that
$$
I_2(t\wedge\tau_R)\leq C_{T}\int^{t\wedge\tau_R}_0\|\w(s)\|^2_{\mH^0}\dif s
$$
and by {\bf (H3)},
$$
I_4(t\wedge\tau_R)\leq \sup_{t\geq 0, x\in\mD}\|\sigma(t,x)\|_{l^2}
\cdot\int^{t\wedge\tau_R}_0\|\nabla\w(s)\|^2_{\mH^0}\dif s
+C_T\int^{t\wedge\tau_R}_0\|\w(s)\|^2_{\mH^0}\dif s.
$$
Taking expectations for (\ref{Lp7}) and combining
the above calculations as well as (\ref{Sigma}), we find that for any $t\in[0,T]$,
\begin{align*}
\mE\|\w(t\wedge\tau_R)\|^2_{\mH^0}\leq C_{R,T}\cdot\mE\left(\int^{t\wedge\tau_R}_0
\|\w(s)\|^2_{\mH^0}\dif s\right)\leq C_{R,T}\int^t_0\mE\|\w(s\wedge\tau_R)\|^2_{\mH^0}\dif s.
\end{align*}
By Gronwall's inequality, we get  for any $t\in[0,T]$,
\begin{align*}
\mE\|\w(t\wedge\tau_R)\|^2_{\mH^0}=0.
\end{align*}
Now the uniqueness follows by letting $R\uparrow\infty$ and Fatou's lemma.
\end{proof}

\subsection{Existence of Martingale Solutions}

We now prove the existence of a weak solution to Eq. (\ref{NS}).
\bt\label{Th1}
Under {{\bf (H1)}}-{\bf (H3)}, for any initial law $\vartheta\in\cP(\mH^1)$,
there exists a weak solution for Eq. (\ref{NS}) in the sense of
Definition \ref{def}.
\et

We shall use Galerkin's approximation to prove this theorem.
In the following, we  fix a stochastic basis $(\Omega,\cF,P; (\cF_t)_{t\geq 0})$,
and an infinite sequence of independent standard $(\cF_t)$-Brownian
motions $\{W^k(t), t\geq 0, k\in\mN\}$, as well as an $\cF_0$-measurable
random variable $\u_0$ having law $\vartheta$.

Recall that $\sE=\{\e_i,i\in\mN\}\subset\cV$ is a complete orthonormal basis  of $\mH^1$.
Set
$$
\mH^1_n:=\mathrm{span}\{e_i, i=1,\cdots,n\}
$$
and for $\u\in\mH^0$,
$$
\Pi_n\u:=\sum_{i=1}^n\<\u,\e_i\>_{\mH^1}\e_i=\sum_{i=1}^n\<\u,(I-\Delta)\e_i\>_{\mH^0}\e_i.
$$
Consider the following finite dimensional stochastic ordinary differential equation
in $\mH^1_n$
\begin{align*}
\left\{
\begin{array}{lcl}
\dif \u_n(t)&=& [\Pi_n A(\u_n(t))+\Pi_n\f(t,\u_n(t))]\dif t+\sum_k \Pi_nB_k(t,\u_n(t))\dif W^k_t,\\
\u_n(0)&=& \Pi_n\u_0.
\end{array}
\right.
\end{align*}
By Lemmas \ref{Le1} and \ref{Le6}, we have, for some $C_{n,N}>0$ and any $\u\in\mH^1_n$,
\begin{align*}
\<\u,\Pi_nA(\u)+\Pi_n\f(t,\u)\>_{\mH^1_n}&\leq C_{n,N}(\|\u\|^2_{\mH^1_n}+1),\\
\|\Pi_nB(t,\u)\|_{l^2\otimes\mH^1_n}&\leq C_{n,N}(\|\u\|^2_{\mH^1_n}+1).
\end{align*}
Moreover, by {\bf (H1)}-{\bf (H3)} it is easy to see that
$$
\mH^1_n\ni\u\mapsto \Pi_nA(\u)+\Pi_n\f(t,\u)\in\mH^1_n
$$
and
$$
\mH^1_n\ni \u\mapsto \Pi_n B (t,\u)\in l^2\times\mH^1_n
$$
are locally Lipschitz continuous. Hence, by the theory of SDE  (cf. \cite{Kr0, Roe}),
there is a unique continuous ($\cF_t$)-adapted process $\u_n(t)$ satisfying
\begin{align}
\u_n(t)&= \u_n(0)+\int^t_0\Pi_nA(\u_n(s))\dif s+\int^t_0\Pi_n\f(s,\u_n(s))\dif s
+\sum_{k=1}^\infty\int^t_0\Pi_nB_k(s,\u_n(s))\dif W^k_s\label{OP4}
\end{align}
and for any $n\geq i$,
\begin{align}
\<\u_n(t),\e_i\>_{\mH^1}&= \<\u_0,\e_i\>_{\mH^1}+\int^t_0\lb A(\u_n(s)),\e_i\rb\dif s
+\int^t_0\<\f(s,\u_n(s)),\e_i\>_{\mH^1}\dif s\no\\
&\quad +\sum_{k=1}^\infty\int^t_0\<B_k(s,\u_n(s)),\e_i\>_{\mH^1}\dif W^k_s.\label{Es6}
\end{align}

We now prove a series of lemmas.
\bl
For any $T>0$, there exists a positive constant $C_{T,N}>0$ such that for any $n\in\mN$,
\begin{align}
\mE\left(\sup_{t\in[0,T]}\|\u_n(t)\|^2_{\mH^1}\right)+\int^T_0\mE\|\u_n(s)\|^2_{\mH^2}\dif s
+\int^T_0\mE\|\nabla|\u_n(s)|^2\|^2_{L^2}\dif s
\leq C_{T,N},\label{Es1}
\end{align}
and also in the periodic case
\begin{align}
\int^T_0\mE\|\u_n(s)\|^4_{L^4}\dif s\leq C_{T,N}.\label{Es11}
\end{align}
\el
\begin{proof}
By It\^o's formula and Lemmas \ref{Le1} and \ref{Le6}, we have
\begin{align}
\|\u_n(t)\|^2_{\mH^1}&= \|\u_0\|^2_{\mH^1}+2\int^t_0\lb A(\u_n(s)),\u_n(s)\rb\dif s
+2\int^t_0\<\f(s,\u_n(s)),\u_n(s)\>_{\mH^1}\dif s\no\\
&\quad +M(t)+\int^t_0\|B(s,\u_n(s))\|_{L_2(l^2;\mH^1)}^2\dif s\no\\
&\leq \|\u_0\|^2_{\mH^1}-\int^t_0\|\u_n(s)\|^2_{\mH^2}\dif s
-\int^t_0\||\u_n(s)|\cdot|\nabla\u_n(s)|\|^2_{L^2}\dif s\label{Es5}\\
&\quad +C\cdot N\int^t_0\|\nabla\u_n(s)\|^2_{\mH^0}\dif s+2\int^t_0\|\u_n(s)\|^2_{\mH^0}\dif s\no\\
&\quad +2\int^t_0\|\f(s,\u_n(s))\|_{\mH^0}\cdot\|\u_n(s)\|_{\mH^2}\dif s+M(t)\no\\
&\quad +\int^t_0\Big(\frac{1}{2}\|\u_n(s)\|^2_{\mH^2}
+C_T\|\u_n(s)\|^2_{\mH^1}+C\|H_{\h}(s)\|_{L^1(\mD)}\Big)\dif s,\no
\end{align}
where $M(t)$ is a continuous martingale defined by
\begin{align*}
M(t):=2\sum_{k=1}^\infty\int^t_0\<B_k(s,\u_n(s)),\u_n(s)\>_{\mH^1}\dif W^k_s.
\end{align*}
Taking expectations and by Young's inequality, one finds that for any $t\in[0,T]$,
\begin{align*}
\mE\|\u_n(t)\|^2_{\mH^1}&\leq \mE\|\u_0\|^2_{\mH^1}-\frac{1}{4}\int^t_0\mE\|\u_n(s)\|^2_{\mH^2}\dif s
-\int^t_0\mE\||\u_n(s)|\cdot|\nabla\u_n(s)|\|^2_{L^2}\dif s\\
&\quad +C\cdot N\int^t_0\mE\|\nabla\u_n(s)\|^2_{\mH^0}\dif s+C_T \int^t_0\mE\|\u_n(s)\|^2_{\mH^0}\dif s\\
&\quad +C_T\int^t_0\Big(\|H_{\f}(s)\|_{L^1(\mD)}+\|H_{\h}(s)\|_{L^1(\mD)}\Big)\dif s.
\end{align*}
Hence, by Gronwall's inequality, we have for any $T>0$,
\begin{align}
\sup_{t\in[0,T]}\mE\|\u_n(t)\|^2_{\mH^1}+\int^T_0\mE\|\u_n(s)\|^2_{\mH^2}\dif s
+\int^T_0\mE\|\nabla|\u_n(s)|^2\|^2_{L^2}\dif s
\leq C_{T,N}.\label{Es21}
\end{align}
Here, the constant $C_{T,N}$ is independent of $n$, and
we have used that $|\nabla|\u|^2|\leq C|\u|\cdot|\nabla\u|$.

Furthermore, from (\ref{Es5}) and using Burkholder's inequality, Young's inequality,
Lemma \ref{Le6} and (\ref{Es21}), we have for any $T>0$ and $\epsilon>0$,
\begin{align*}
\mE\left(\sup_{t\in[0,T]}\|\u_n(t)\|^2_{\mH^1}\right)&\leq C_{T,N}+
C\mE\left(\int^T_0\|B(s,\u_n(s))\|^2_{L_2(l^2;\mH^0)}\cdot\|\u_n(s)\|^2_{\mH^2}\dif s\right)^{1/2}\no\\
&\leq C_{T,N}+\epsilon\cdot\mE\left(\sup_{t\in[0,T]}\|B(s,\u_n(s))\|^2_{L_2(l^2;\mH^0)}\right)
+C_\epsilon\int^T_0\mE\|\u_n(s)\|^2_{\mH^2}\dif s\no\\
&\leq C_{T,N,\epsilon}+\epsilon\cdot C_T\mE\left(\sup_{t\in[0,T]}\|\u_n(t)\|^2_{\mH^1}\right).
\end{align*}
Choosing $\epsilon$  small enough, we get
\begin{align*}
\mE\left(\sup_{t\in[0,T]}\|\u_n(t)\|^2_{\mH^1}\right)\leq C_{T,N}.
\end{align*}

In the periodic case, since $\sE$ is also orthogonal in $\mH^0$, we have by (\ref{Es44}) and
{\bf (H1)}-{\bf (H3)},
\begin{align*}
\mE\|\u_n(t)\|^2_{\mH^0}&= \mE\|\u_0\|^2_{\mH^0}+2\int^t_0\mE\< A(\u_n(s)),\u_n(s)\>_{\mH^0}\dif s\\
&\quad +2\int^t_0\mE\<\f(s,\u_n(s)),\u_n(s)\>_{\mH^0}\dif s
+\int^t_0\mE\|B(s,\u_n(s))\|_{L_2(l^2;\mH^0)}^2\dif s\no\\
&\leq \mE\|\u_0\|^2_{\mH^1}-2\int^t_0\mE\|\nabla\u_n(s)\|^2_{\mH^0}\dif s
-2\int^t_0\mE\|\u_n(s)\|^4_{L^4}\dif s\\
&\quad +C_N\int^t_0\mE\|\u_n(s)\|^2_{\mH^0}\dif s+C_T,
\end{align*}
which yields (\ref{Es11}) by Gronwall's lemma.
\end{proof}
\bl
Let $\mu_n$ be the law of $\u_n$ in $(\mX,\cB(\mX))$. Then the family of probability measures
$\{\mu_n,n\in\mN\}$ is tight on $(\mX,\cB(\mX))$.
\el
\begin{proof}
Set for $R>0$
$$
\tau^n_R:=\inf\{t\geq 0: \|\u_n(t)\|_{\mH^1}\geq R\}.
$$
Then, by (\ref{Es1}) we have for any $T>0$,
\begin{align}
\sup_n P(\tau^n_R<T)=\sup_n P\left(\sup_{t\in[0,T]}\|\u_n(t)\|_{\mH^1}>R\right)\leq \frac{C_{T,N}}{R^2}.\label{Es88}
\end{align}
On the other hand, from (\ref{Es6}) and using (\ref{Op1}), Lemma \ref{Le6} and Burkholder's inequality,
we have for any $q\geq 2$ and $s,t\in[0,T]$, $\e\in\sE$,
\begin{align*}
&\quad \mE|\<\u_n(t\wedge\tau^n_R)-\u_n(s\wedge\tau^n_R),\e\>_{\mH^1}|^q\\
&\leq C\mE\left|\int^{t\wedge\tau^n_R}_{s\wedge\tau^n_R}\lb A(\u_n(s)),\e\rb\dif s\right|^q
+C\mE\left|\int^{t\wedge\tau^n_R}_{s\wedge\tau^n_R}\<\f(s,\u_n(s)),\e\>_{\mH^1}\dif s\right|^q\no\\
&\quad +C\mE\left|\int^{t\wedge\tau^n_R}_{s\wedge\tau^n_R}\<B_k(s,\u_n(s)),\e\>_{\mH^1}\dif W^k_s\right|^q\\
&\leq C_{\e}\cdot\mE\left|\int^{t\wedge\tau^n_R}_{s\wedge\tau^n_R}(1+\|\u_n(s)\|^3_{\mH^1})\dif s\right|^q
+C_{\e}\cdot\mE\left|\int^{t\wedge\tau^n_R}_{s\wedge\tau^n_R}\|\f(s,\u_n(s))\|_{\mH^0}\dif s\right|^q\no\\
&\quad +C_{\e}\cdot\mE\left|\int^{t\wedge\tau^n_R}_{s\wedge\tau^n_R}
\|B(s,\u_n(s))\|^2_{L_2(l^2;\mH^0)}\dif s\right|^{q/2}\leq C_{\e,R,T}\cdot|t-s|^{q/2}.
\end{align*}
By Kolomogorov's criterion (cf. \cite{Ka}), we get for any $T>0$ and $0<\a<\frac{1}{2}$,
\begin{align*}
\mE\left(\sup_{ s,t\in[0,T],|t-s|\leq \delta}
|\<\u_n(t\wedge\tau^n_R)-\u_n(s\wedge\tau^n_R),\e\>_{\mH^1}|\right)\leq C_{\e,R,T}\cdot\delta^{\a}.
\end{align*}
So, for any $\epsilon>0$ and $R>0$,
\begin{align*}
&\sup_n P\left\{\sup_{s,t\in[0,T],|t-s|\leq \delta}
|\<\u_n(t)-\u_n(s),\e\>_{\mH^1}|>\epsilon\right\}\\
&\qquad\leq \sup_n P\left\{\sup_{s,t\in[0,T],|t-s|\leq \delta}
|\<\u_n(t)-\u_n(s),\e\>_{\mH^1}|>\epsilon;\tau^n_R\geq T\right\}+\sup_n P\left\{\tau^n_R<T\right\}\\
&\qquad\leq\frac{C_{\e,R,T}\cdot\delta}{\epsilon}+\frac{C_{T,N}}{R^2},
\end{align*}
which then gives that
\begin{align}
\lim_{\delta\downarrow 0}\sup_n P\left\{\sup_{s,t\in[0,T],|t-s|\leq \delta}
|\<\u_n(t)-\u_n(s),\e\>_{\mH^1}|>\epsilon\right\}=0.\label{Es9}
\end{align}
The tightness of $\{\mu_n,n\in\mN\}$ now follows from
(\ref{Es88}), (\ref{Es9}) and Lemma \ref{Le5}.
\end{proof}

In the sequel, without loss of generality, we assume that $\mu_n$
weakly converges to $\mu\in\cP(\mX)$.
By Skorohod's embedding theorem (cf. \cite{Ka}), there exist a probability space $(\tilde\Omega,\tilde\cF,\tilde P)$
and $\mX$-valued random variables $\tilde \u^n$ and $\tilde \u$ such that

(I) $\tilde \u^n$ has the same law as $\u^n$ in $\mX$ for each $n\in\mN$;

(II) $\tilde \u^n\rightarrow \tilde \u$ in $\mX$, $\tilde P$-a.e., and $\tilde \u$ has law $\mu$.

Moreover, by (\ref{Es1}) and Fatou's lemma, we have  for any $T>0$,
\begin{align}
&\quad \mE^{\mu}\left(\sup_{t\in[0,T]}\|\u(t)\|^2_{\mH^1}\right)=
\mE^{\tilde P}\left(\sup_{t\in[0,T]}\|\tilde \u(t)\|^2_{\mH^1}\right)<+\infty,\label{Es009}\\
&\quad \int^T_0\mE^{\mu}\|\u(s)\|^2_{\mH^2}\dif s=\int^T_0
\mE^{\tilde P}\|\tilde \u(s)\|^2_{\mH^2}\dif s<+\infty,\\
&\quad \int^T_0\mE^{\mu}\|\nabla|\u(s)|^2\|^2_{L^2}\dif s
=\int^T_0\mE^{\tilde P}\|\nabla|\tilde \u(s)|^2\|^2_{L^2}\dif s< +\infty\label{Es09}
\end{align}
and also in the periodic case
\begin{align}
\int^T_0\mE^{\mu}\|\u(s)\|^4_{L^4}\dif s
=\int^T_0\mE^{\tilde P}\|\tilde \u(s)\|^4_{L^4}\dif s< +\infty.\label{Es10}
\end{align}

Let  $h\in C^\infty_0(\mR)$ and $\e\in\sE$. Define for any $t\geq 0$ and $\u\in\mX$,
$$
M^h_{\e}(t,\u):=I^h_1(t,\u)-I^h_2(t,\u)-I^h_3(t,\u)-I^h_4(t,\u)-I^h_5(t,\u),
$$
where
\begin{align*}
I^h_1(t,\u)&:= h(\<\u(t),\e\>_{\mH^1}),\\
I^h_2(t,\u)&:= h(\<\u(0),\e\>_{\mH^1}),\\
I^h_3(t,\u)&:= \int^t_0h'(\<\u(s),\e\>_{\mH^1})\cdot\lb  A(\u(s)),\e\rb\dif s,\\
I^h_4(t,\u)&:= \int^t_0h'(\<\u(s),\e\>_{\mH^1})\cdot \<\f(s,\u(s)),\e\>_{\mH^1}\dif s,\\
I^h_5(t,\u)&:= \frac{1}{2}\int^t_0h''(\<\u(s),\e\>_{\mH^1})
\cdot\|\<B(s,\u(s)),\e\>_{\mH^1}\|^2_{l^2}\dif s.
\end{align*}

Note that $\e\in\sE\subset\cV$ has compact support, there exists $m\in\mN$ such that
\begin{align}
\mathrm{supp}\{\e\}\subset\cO:=\{x\in\mR^3, |x|\leq m\}.\label{Sup}
\end{align}
\bl\label{Le7}
 We have
\begin{align}
\sup_n\mE^{\tilde P}|M^h_{\e}(t,\tilde\u_n)|^{4/3}
+\mE^{\tilde P}|M^h_{\e}(t,\tilde\u)|^{4/3}<+\infty.\label{Op9}
\end{align}
\el
\begin{proof}
It is clear that $I^h_1(t,\tilde\u_n)$ and $I^h_2(t,\tilde\u_n)$ are bounded by some constant $C_h$.
For $I^h_3$,  noting that in the whole space case,
$$
\|\u\|_{L^6}\leq C\|\nabla\u\|_{L^2},
$$
 by (\ref{Op1}) and (\ref{Es1}), we have
\begin{align*}
\mE^{\tilde P}|I^h_3(t,\tilde\u_n)|^{4/3}
&\leq C_{T,h}\int^T_0\mE^{\tilde P}\big|\lb  A(\tilde\u_n(s)),\e\rb\big|^{4/3}\dif s\\
&\leq C_{T,h,\e}\int^T_0\mE^{\tilde P}\left(1+\|\tilde\u_n(s)\|^4_{L^3(\cO)}\right)\dif s\\
&\leq C_{T,h,\e}\int^T_0\mE^{\tilde P}\left(1+\|\tilde\u_n(s)\|^4_{L^{12}(\cO)}\right)\dif s\\
&= C_{T,h,\e}\int^T_0\mE^{\tilde P}\left(1+\||\tilde\u_n(s)|^2\|^2_{L^6}\right)\dif s\\
&\leq C_{T,h,\e}\int^T_0\mE^{\tilde P}\left(1+\|\nabla|\tilde\u_n(s)|^2\|^2_{L^2}\right)\dif s\\
&\leq C_{T,h,\e,N}.
\end{align*}
In the periodic case, by (\ref{Op1}) and (\ref{Es11}), we have
\begin{align*}
\mE^{\tilde P}|I^h_3(t,\tilde\u_n)|^{4/3}
&\leq C_{T,h}\int^T_0\mE^{\tilde P}\left(1+\|\tilde\u_n(s)\|^4_{L^3(\mT^3)}\right)\dif s\\
&\leq C_{T,h,\e}\int^T_0\mE^{\tilde P}\left(1+\|\tilde\u_n(s)\|^4_{L^4(\mT^3)}\right)\dif s\\
&\leq C_{T,h,\e,N}.
\end{align*}
For $I^h_5$, by (\ref{Op2}), we similarly have
$$
\mE^{\tilde P}|I^h_5(t,\tilde\u_n)|^2\leq C_{T,h,\e,N}.
$$
For  $I^h_4$, it is clear that
$$
\mE^{\tilde P}|I^h_4(t,\tilde\u_n)|^2\leq C_{T,h,\e,N}.
$$
Moreover, by (\ref{Es009})-(\ref{Es10}), we also have
$$
\mE^{\tilde P}|M^h_{\e}(t,\tilde\u)|^{4/3}<+\infty.
$$
The proof is thus complete.
\end{proof}
\bl
For any $t>0$ and $\epsilon>0$,
\begin{align}
\lim_{n\rightarrow\infty}\tilde P\big(|M^h_{\e}(t,\tilde\u_n)
-M^h_{\e}(t,\tilde\u)|>\epsilon\big)=0.\label{Es99}
\end{align}
That is, $M^h_{\e}(t,\tilde\u_n)$ converges to $M^h_{\e}(t,\tilde\u)$ in probability $\tilde P$
as $n\rightarrow\infty$.
\el
\begin{proof}
Recalling the definition of $\mX$ in Subsection 2.2, by (II) we have
$$
\lim_{n\to\infty}\int_{\cO}|\tilde \u_n(t,x,\tilde\omega)-\tilde \u(t,x,\tilde\omega)|^2\dif x=0,\ \
\tilde P-a.a.\ \tilde\omega\in\tilde\Omega,
$$
where $\cO$ is from (\ref{Sup}).
Thus, by the dominated convergence theorem, we have
\begin{align*}
&\quad \lim_{n\rightarrow\infty}\mE^{\tilde P}|I^h_1(t,\tilde\u_n)-I^h_1(t,\tilde\u)|=0,\\
&\quad \lim_{n\rightarrow\infty}\mE^{\tilde P}|I^h_2(t,\tilde\u_n)-I^h_2(t,\tilde\u)|=0.
\end{align*}
For$I^h_3$, define for any $R>0$,
$$
\tilde\tau^n_R:=\inf\{t\geq 0: \|\tilde\u_n(t)\|_{\mH^1}\geq R\}.
$$
Then, by (I) and (\ref{Es1}), for any $T>0$, we have
$$
\sup_n \tilde P(\tilde\tau^n_R\leq T)\leq \frac{C_{T,N}}{R^2}.\label{Es8}
$$
Thus, by the dominated convergence theorem and Lemma \ref{Le2}, we have from the proof of
Lemma \ref{Le7},
\begin{align*}
&\lim_{n\rightarrow\infty}\tilde P(|I^h_3(t,\tilde\u_n)-I^h_3(t,\tilde\u)|>\epsilon)\\
&\quad \leq \lim_{R\rightarrow\infty}\lim_{n\rightarrow\infty}\tilde
P(|I^h_3(t,\tilde\u_n)-I^h_3(t,\tilde\u)|>\epsilon; \tilde\tau^n_R>t)
+\lim_{R\rightarrow\infty}\sup_n \tilde P(\tilde\tau^n_R\leq T)\\
&\quad \leq \lim_{R\rightarrow\infty}\lim_{n\rightarrow\infty}
\mE^{\tilde P}\left(1_{\{\tilde\tau^n_R>t\}}\cdot|I^h_3(t,\tilde\u_n)-I^h_3(t,\tilde\u)|\right)/\epsilon\\
&\quad \leq \lim_{R\rightarrow\infty}\mE^{\tilde P}\Bigg(\int^t_0\varlimsup_{n\rightarrow\infty}
\Big(1_{\{\tilde\tau^n_R>t\}}\cdot\Big|h'(\<\tilde\u_n(s),\e\>_{\mH^1})\cdot\lb  A(\tilde\u_n(s)),\e\rb\\
&\quad \quad -h'(\<\tilde\u(s),\e\>_{\mH^1})\cdot\lb  A(\tilde\u(s)),\e\rb\Big|\Big)\dif s\Bigg)/\epsilon=0.
\end{align*}
Similarly, we also have
\begin{align*}
&\quad \lim_{n\rightarrow\infty}\tilde P(|I^h_4(t,\tilde\u_n)-I^h_4(t,\tilde\u)|>\epsilon)=0,\\
&\quad \lim_{n\rightarrow\infty}\tilde P(|I^h_5(t,\tilde\u_n)-I^h_5(t,\tilde\u)|>\epsilon)=0.
\end{align*}
Combining the above calculations yields (\ref{Es99}).
\end{proof}

We can now give the proof of Theorem \ref{Th1}.

\vspace{4mm}

{\it Proof of Theorem \ref{Th1}}:
Now let $t>s$ and $G$ be any bounded and real valued
$\cB_s(\mX)$-measurable continuous function on $\mX$. Then by (\ref{Op9}) and  (\ref{Es99}), we have
\begin{align*}
\mE^\mu\Big((M^h_{\e}(t,\u)-M^h_{\e}(s,\u))\cdot G(\u)\Big)
&= \mE^{\tilde P}\Big((M^h_{\e}(t,\tilde\u)-M^h_{\e}(s,\tilde\u))\cdot G(\tilde\u)\Big)\\
&= \lim_{n\rightarrow\infty}
\mE^{\tilde P}\Big((M^h_{\e}(t,\tilde\u_n)-M^h_{\e}(s,\tilde\u_n))\cdot G(\tilde\u_n)\Big)\\
&= \lim_{n\rightarrow\infty}
\mE^{P}\Big((M^h_{\e}(t,\u_n)-M^h_{\e}(s,\u_n))\cdot G(\u_n)\Big)=0,
\end{align*}
where the last step is due to the martingale property of $M^h_{\e}(t,\u_n)$
on $(\Omega,\cF,P; (\cF_t)_{t\geq 0})$ and $G(\u_n)\in\cF_s$.
This means that
$\{M^h_{\e}(t,\u),t\geq0\}$ is a $\cB_t(\mX)$-martingale.
The existence of a weak solution to Eq. (\ref{NS}) now follows from Proposition \ref{Th2}.

\vspace{5mm}

Summarizing Theorems \ref{Uni}, \ref{Th1} and \ref{Ya-Wa}, we have
the following main result in the present paper.
\bt\label{main}
Under {{\bf (H1)}}-{\bf (H3)}, for any $\u_0\in\mH^1$,
there exists a unique $\u(t,x)$
such that
\begin{enumerate}[($1^o$)]
\item $\u\in L^2(\Omega,P; C([0,T],\mH^1))\cap L^2(\Omega,P; L^2([0,T],\mH^2))$ for any $T>0$, and
\begin{align}
\mE\left(\sup_{t\in[0,T]}\|\u(t)\|^2_{\mH^1}\right)+\int^T_0\mE\|\u(s)\|^2_{\mH^2}\dif s
\leq C_{T}(1+\|\u_0\|^2_{\mH^1})N;\label{Op10}
\end{align}
\item it holds that in $\mH^0$,
\begin{align*}
\u(t)=\u_0+\int^t_0\Big[A(\u(s))+\sP\f(s,\u(s))\Big]\dif s
+\sum_{k=1}^\infty\int^t_0 B_k(s,\u(s))\dif W^k_s,
\end{align*}
for all $t\geq 0$, $P$-a.s..
\end{enumerate}
\et
\begin{proof}
We only need to prove estimate (\ref{Op10}). By It\^o's formula, (\ref{Es444}) and Lemma \ref{Le6}, we have
\begin{align*}
\mE\|\u(t)\|^2_{\mH^0}&= \|\u_0\|^2_{\mH^0}+2\int^t_0\mE\<A(\u(s)),\u(s)\>_{\mH^0}\dif s\\
&\quad +2\int^t_0\mE\<\f(s,\u(s)),\u(s)\>_{\mH^0}\dif s+\int^t_0\mE\|B(s,\u(s))\|^2_{L_2(l^2;\mH^0)}\dif s\\
&\leq \|\u_0\|^2_{\mH^0}+C-\frac{1}{2}\int^t_0\mE\|\u(s)\|^2_{\mH^1}\dif s
+C\int^t_0\mE\|\u(s)\|^2_{\mH^0}\dif s.
\end{align*}
By Gronwall's inequality, we obtain
\begin{align*}
\sup_{t\in[0,T]}\mE\|\u(t)\|^2_{\mH^0}+\int^T_0\mE\|\u(s)\|^2_{\mH^1}\dif s\leq C_T(\|\u_0\|^2_{\mH^0}+1).
\end{align*}
Using this estimate, as in the proof of  (\ref{Es1}), we obtain (\ref{Op10}).
\end{proof}

\section{Feller Properties and  Invariant Measures}

In the following, we consider the time homogenous case, i.e., the coefficients $\f$, $\sigma$ and $\h$
are independent of $t$, and assume a stronger assumption than {\bf (H3)}, namely:
\begin{enumerate}[{\bf (H3)$'$}]
\item There exist a constant $C_{\h}>0$ and a
 function $H_{\h}(x)\in L^1(\mD)$ such that for any $x\in\mD,\u,\v\in\mR^3$
and $j=1,2,3$,
\begin{align*}
\|\p_{x^j}\h(x,\u)\|_{l^2}^2+\|\h(x,\u)\|_{l^2}^2&\leq C_{\h}\cdot|\u|^2+H_{\h}(x),\\
\|\p_{x^j}\h(x,\u)-\p_{x^j}\h(x,\v)\|_{l^2}&\leq C_{\h}\cdot|\u-\v|,\\
\|\p_{u^j}\h(x,\u)\|_{l^2}&\leq C_{\h},\\
\|\p_{u^j}\h(x,\u)-\p_{v^j}\h(x,\v)\|_{l^2}&\leq C_{\h}\cdot|\u-\v|.
\end{align*}
\end{enumerate}

For fixed initial value $\u_0=\v\in\mH^1$, we denote the unique solution in
Theorem \ref{main} by $\u(t;\v)$. Then $\{\u(t;\v): \v\in\mH^1, t\geq 0\}$ forms a strong Markov process
with state space $\mH^1$. We have:
\bl\label{Le8}
For $\v,\v'\in\mH^1$ and $R>0$, define
$$
\tau^{\v}_R:=\inf\left\{t\geq 0: \|\u(t;\v)\|_{\mH^1}>R\right\}
$$
and
$$
\tau^{\v,\v'}_R:=\tau^{\v}_R\wedge\tau^{\v'}_R.
$$
Assume {\bf (H1)}, {\bf (H2)} and {\bf (H3)$'$}, then
\begin{align*}
\mE\|\u(t\wedge\tau^{\v,\v'}_R;\v)-\u(t\wedge\tau^{\v,\v'}_R;\v')\|^2_{\mH^1}
\leq C_{t,R}\cdot\|\v-\v'\|_{\mH^1}^2.
\end{align*}
\el
\begin{proof}
Write $\u(t):=\u(t;\v)$, $\tilde\u(t):=\u(t,\v')$ and
$$
\w(t):=\u(t)-\tilde\u(t).
$$
Set $t_R:=\tau^{\v,\v'}_R\wedge t$.
By It\^o's formula (cf. \cite{Ro,Roe}), we have
\begin{align*}
\|\w(t_R)\|^2_{\mH^1}&= \|\w(0)\|_{\mH^1}^2+2\int^{t_R}_0\<A(\u(s))-A(\tilde\u(s)),\w(s)\>_{\mH^1}\dif s\\
&\quad +2\int^{t_R}_0\<\f(s,\u(s))-\f(s,\tilde\u(s)),\w(s)\>_{\mH^1}\dif s\\
&\quad +2\sum_{k=1}^\infty\int^{t_R}_0\<B_k(s,\u(s))-B_k(s,\tilde\u(s)),\w(s)\>_{\mH^1}\dif W^k_s\\
&\quad +\sum_{k=1}^\infty\int^{t_R}_0\|B_k(s,\u(s))-B_k(s,\tilde\u(s))\|^2_{\mH^1}\dif s\\
&=:\|\w(0)\|_{\mH^1}^2+I_1(t_R)+I_2(t_R)+I_3(t_R)+I_4(t_R).
\end{align*}
By $|g_N(r)-g_N(r')|\leq |r-r'|$ and Young's inequality, it is easy to see that
\begin{align*}
I_1(t_R)&= -2\int^{t_R}_0\|\w(s)\|^2_{\mH^2}\dif s+2\int^{t_R}_0\|\w(s)\|^2_{\mH^1}\dif s\\
&\quad +2\int^{t_R}_0\<((\u(s)\cdot\nabla)\u(s)-(\tilde\u(s)\cdot\nabla)\tilde\u(s)),(I-\Delta)\w(s)\>_{\mH^0}\dif s\\
&\quad -2\int^{t_R}_0\<g_N(|\u(s)|^2)\u(s)-g_N(|\tilde\u(s)|^2)\tilde\u(s),(I-\Delta)\w(s)\>_{\mH^0}\dif s\\
&\leq -\int^{t_R}_0\|\w(s)\|^2_{\mH^2}\dif s+2\int^{t_R}_0\|\w(s)\|^2_{\mH^1}\dif s\\
&\quad +C\int^{t_R}_0\|(\w(s)\cdot\nabla)\u(s)\|^2_{L^2}\dif s
+C\int^{t_R}_0\|(\tilde\u(s)\cdot\nabla)\w(s)\|^2_{L^2}\dif s\\
&\quad +C\int^{t_R}_0\||\w(s)|\cdot(|\u(s)|^2+|\tilde\u(s)|^2)\|^2_{L^2}\dif s.
\end{align*}
By H\"older's inequality and the Sobolev inequality (\ref{Sob}), we further have
\begin{align*}
I_1(t_R)&\leq -\int^{t_R}_0\|\w(s)\|^2_{\mH^2}\dif s+2\int^{t_R}_0\|\w(s)\|^2_{\mH^1}\dif s\\
&\quad +C_R\int^{t_R}_0\|\w(s)\|^2_{L^\infty}\dif s
+C\int^{t_R}_0\|\tilde\u(s)\|^2_{L^6}\cdot\|\nabla\w(s)\|^2_{L^3}\dif s\\
&\quad +C\int^{t_R}_0\|\w(s)\|^2_{L^6}\cdot(\|\u(s)\|_{L^6}^2+\|\tilde\u(s)\|_{L^6}^2)^2\dif s\\
&\leq -\int^{t_R}_0\|\w(s)\|^2_{\mH^2}\dif s+C_R\int^{t_R}_0\|\w(s)\|^2_{\mH^1}\dif s\\
&\quad +C_R\int^{t_R}_0\|\w(s)\|^{3/2}_{\mH^2}\cdot\|\w(s)\|^{1/2}_{\mH^0}\dif s
+C_R\int^{t_R}_0\|\w(s)\|_{\mH^2}\cdot\|\w(s)\|_{\mH^1}\dif s\\
&\leq -\frac{3}{4}\int^{t_R}_0\|\w(s)\|^2_{\mH^2}\dif s+C_R\int^{t_R}_0\|\w(s)\|^2_{\mH^1}\dif s.
\end{align*}
By {\bf (H1)}, {\bf (H2)} and  {\bf (H3)$'$}, we similarly have
\begin{align*}
I_2(t_R)&\leq \frac{1}{4}\int^{t_R}_0\|\w(s)\|^2_{\mH^2}\dif s
+C_R\int^{t_R}_0\|\w(s)\|^2_{\mH^0}\dif s,\\
I_4(t_R)&\leq \frac{1}{2}\int^{t_R}_0\|\w(s)\|^2_{\mH^2}\dif s
+C_R\int^{t_R}_0\|\w(s)\|^2_{\mH^1}\dif s.
\end{align*}

So,
\begin{align*}
\mE\|\w(t\wedge\tau_R)\|^2_{\mH^1}&\leq \|\w(0)\|_{\mH^1}^2+C_R\int^{t_R}_0\|\w(s)\|^2_{\mH^1}\dif s\\
&\leq \|\v-\v'\|_{\mH^1}^2+C_R\int^t_0\|\w(s\wedge\tau_R)\|^2_{\mH^1}\dif s.
\end{align*}
By Gronwall's inequality, we get the desired estimate.
\end{proof}

Let $C_b^{loc}(\mH^1)$ denote the set of all  bounded and
locally uniformly continuous functions on $\mH^1$. Then $C_b^{loc}(\mH^1)$ is clearly a Banach space
under the sup  norm
$$
\|\phi\|_\infty:=\sup_{\u\in\mH^1}|\phi(\u)|.
$$
For $t>0$, we define the semigroup $\bT_t$ associated with
$\{\u(t;\v): \v\in\mH^1, t\geq 0\}$ by
\begin{align*}
\bT_t\phi(\v):=\mE(\phi(\u(t;\v))), \quad \phi\in C_b^{loc}(\mH^1).
\end{align*}

We have:
\bt
Under {\bf (H1)}, {\bf (H2)} and {\bf (H3)$'$}, for every $t>0$,
$\bT_t$ maps $C_b^{loc}(\mH^1)$ into $C_b^{loc}(\mH^1)$.
That is, $(\bT_t)_{t\geq 0}$ is a Feller semigroup on $C_b^{loc}(\mH^1)$.
\et
\begin{proof}
Let $\phi\in C_b^{loc}(\mH^1)$ be given. We want to prove that for any $t>0$ and $m\in\mN$
\begin{align}
\lim_{\delta\rightarrow 0}\sup_{ \v,\v'\in \mB_m, \|\v-\v'\|_{\mH^1}\leq \delta}
|\bT_t\phi(\v)-\bT_t\phi(\v')|=0,\label{PP6}
\end{align}
where $\mB_m:=\{\v\in\mH^1: \|\v\|_{\mH^1}\leq m\}$ denotes the ball in $\mH^1$.

For any $\v,\v'\in\mB_m$ and $R>m$, as in  Lemma \ref{Le8}, define
$$
\tau^{\v}_R:=\{t\geq 0: \|\u(t;\v)\|_{\mH^1}>R\}
$$
and
$$
\tau^{\v,\v'}_R:=\tau^{\v}_R\wedge\tau^{\v'}_R.
$$
By (\ref{Op10}), we have
\begin{align*}
\mE|\phi(\u(t;\v))-\phi(\u(t\wedge\tau^{\v,\v'}_R;\v))|
&\leq 2\|\phi\|_\infty\cdot P(\tau^{\v,\v'}_R<t)
\leq 2\|\phi\|_\infty\cdot  \sup_{\v\in\mB_m}\mE\left(\sup_{s\in[0,t]}\|\u(s;\v)\|^2_{\mH^1}\right)/R^2\\
&\leq 2\|\phi\|_\infty\cdot C_{t,m,N}/R^2.
\end{align*}
For any $\epsilon>0$, choose $R>m$ sufficiently large such that for any $\v,\v'\in\mB_m$
\begin{align}
\mE|\phi(\u(t;\v))-\phi(\u(t\wedge\tau^{\v,\v'}_R;\v))|&\leq \epsilon,\label{PP1}\\
\mE|\phi(\u(t;\v'))-\phi(\u(t\wedge\tau^{\v,\v'}_R;\v'))|&\leq \epsilon.\label{PP2}
\end{align}
For this $R$, since $\phi$ is uniformly continuous on $\mB_R$, one may choose $\eta>0$
such that for any $\u,\u'\in\mB_R$ with $\|\u-\u'\|_{\mH^1}\leq \eta$
\begin{align*}
|\phi(\u)-\phi(\u')|\leq \epsilon.
\end{align*}
Thus, for any $\v,\v'\in\mB_m$ with $\|\v-\v'\|_{\mH^1}\leq \frac{\sqrt{\epsilon}\cdot\eta}
{\sqrt{2C_\phi}\cdot C_{t,R}}$, by Lemma \ref{Le8} we have
\begin{align}
&\mE|\phi(\u(t\wedge\tau^{\v,\v'}_R;\v))-\phi(\u(t\wedge\tau^{\v,\v'}_R;\v'))|\no\\
&\quad \leq \epsilon+2C_\phi\cdot P(\|\u(t\wedge\tau^{\v,\v'}_R;\v)-\u(t\wedge\tau^{\v,\v'}_R;\v')\|_{\mH^1}>\eta)
\leq 2\epsilon.\label{PP3}
\end{align}
Combining (\ref{PP1}) (\ref{PP2}) and (\ref{PP3}), we get (\ref{PP6}).
\end{proof}

In the periodic case, we have the following existence of
invariant measures associated to $(\bT_t)_{t\geq0}$.
\bt\label{Th3}
Under {\bf (H1)}, {\bf (H2)} and {\bf (H3)$'$}, in the periodic
case, there is an invariant measure
$\mu\in\cP(\mH^1)$ associated to the semigroup $(\bT_t)_{t\geq0}$ such that for 
any $t\geq0$ and $\phi\in C_b^{loc}(\mH^1)$,
\begin{align*}
\int_{\mH^1}\bT_t\phi(\u)\mu(\dif \u)=\int_{\mH^1}\phi(\u)\mu(\dif \u).
\end{align*}
\et
\begin{proof}
In the following, we assume that $\u_0=0$.
Using It\^o's formula, we have by (\ref{Es44}) and (\ref{LP1}),
\begin{align*}
\mE\|\u(t)\|^2_{\mH^0}&= 2\int^t_0\mE\<A(\u(s)),\u(s)\>_{\mH^0}\dif s
+2\int^t_0\mE\<\f(\u(s)),\u(s)\>_{\mH^0}\dif s
+\int^t_0\mE\|B(\u(s))\|^2_{L_2(l^2;\mH^0)}\dif s\\
&\leq -\frac{3}{2}\int^t_0\mE\|\u(s)\|^2_{\mH^1}\dif s
-2\int^t_0\mE\|\u(s)\|^4_{L^4}\dif s+C_{\h,\f,N}\int^t_0\mE\|\u(s)\|^2_{\mH^0}\dif s+C_{\h,\f}\cdot t.
\end{align*}
In the periodic case, noting that for any $\epsilon>0$
$$
\|\u\|^2_{\mH^0}\leq C\|\u\|^2_{L^4}\leq \epsilon\|\u\|^4_{L^4}+C_\epsilon,
$$
we further have
\begin{align*}
\mE\|\u(t)\|^2_{\mH^0}\leq-\frac{3}{2}\int^t_0\mE\|\u(s)\|^2_{\mH^1}\dif s
-\int^t_0\mE\|\u(s)\|^4_{L^4}\dif s+C_{\h,\f,N}\cdot t.
\end{align*}
Hence, for any $t\geq 0$
\begin{align}
\mE\|\u(t)\|^2_{\mH^0}+\int^t_0\mE\|\u(s)\|^2_{\mH^1}\dif s
+\int^t_0\mE\|\u(s)\|^4_{L^4}\dif s\leq C_{\h,\f,N}\cdot t.\label{PO3}
\end{align}
On the other hand, by It\^o's formula again and (\ref{Es4}), (\ref{LP2}), as above we have
\begin{align*}
\mE\|\u(t)\|^2_{\mH^1}&= 2\int^t_0\mE\lb A(\u(s)),\u(s)\rb\dif s
+2\int^t_0\mE\<\f(\u(s)),\u(s)\>_{\mH^1}\dif s+\int^t_0\mE\|B(\u(s))\|^2_{L_2(l^2;\mH^1)}\dif s\\
&\leq -\frac{1}{4}\int^t_0\mE\|\u(s)\|^2_{\mH^2}\dif s+
C_{\h,\f,N}\cdot \int^t_0\mE\|\u(s)\|^2_{\mH^1}\dif s+C_{\h,\f,N}\cdot t\\
&\leq -\frac{1}{4}\int^t_0\mE\|\u(s)\|^2_{\mH^2}\dif s+C_{\h,\f,N}\cdot t.
\end{align*}
Therefore, for any $t\geq 0$
\begin{align*}
\frac{1}{t}\int^t_0\mE\|\u(s)\|^2_{\mH^2}\dif s\leq C_{\h,\f,N}.
\end{align*}
In the periodic case, since $\mH^2$ is compactly embedded into $\mH^1$,  the existence of
an invariant measure $\mu$
now follows from the classical Krylov-Bogoliubov method (cf. \cite{DaZa}).
\end{proof}

\section{Ergodicity: Uniqueness of Invariant Measures}

In the following, we shall work in the case of $\mD=\mT^3$, and
suppose that for $\f\in\mH^0$, the mean value of $f$ on $\mT^3$ vanishes, i.e.,
$$
\int_{\mT^3}\f(x)\dif x=0.
$$
In this case, we assume that the orthonormal basis $\sE$
of $\mH^1$ consists of the eigenvectors of $\sP\Delta$, i.e,
$$
\sP\Delta\e_i=-\lambda_i\e_i,\ \ \<\e_i,\e_i\>_{\mH^1}=1,\  \ i=1,2,\cdots,
$$
where $0<\lambda_1\leq\cdots\leq\lambda_n\uparrow\infty$.
Recalling that the following Poincare inequality holds:
\begin{align}
\|\u\|_{\mH^0}^2\leq 1/\lambda_1\|\nabla\u\|_{\mH^0}^2,\label{Poi}
\end{align}
two equivalent norms in $\mH^1$ and $\mH^2$ are given by
$$
\|\u\|_{\mH^1}:=\|\nabla\u\|_{\mH^0},\ \ \|\u\|_{\mH^2}:=\|\Delta\u\|_{\mH^0}.
$$
We shall use these two norms in what follows.

For $m\in\mN$, let $\Omega:=C_0(\mR_+;\mR^m)$ denote the space of all continuous functions with
initial values $0$, $P$ the standard Wiener measure on $\cF:=\cB(C_0(\mR_+;\mR^m))$.
Then, the coordinate process
$$
W_t(\om):=\om(t),\ \ \om\in\Omega,
$$
is a standard Wiener process on $(\Omega,\cF,P)$.

Consider the following stochastic tamed 3D Navier-Stokes equation:
\begin{align}
\label{NSE00}\left\{
\begin{array}{lcl}
\dif \u(t)&=& A(\u(t))\dif t+\dif \w(t),\\
\u(0)&=& \u_0\in\mH^1,
\end{array}
\right.
\end{align}
where $\w(t):=Q W_t$ is the noise,
and the linear map $Q:\mR^m\to \mH^1$ is given by
$$
Qe_i=q_i\e_i,\ \ q_i>0,\ \ i=1,\cdots,m.
$$
Here, $\{e_i,i=1,\cdots,m\}$ is the canonical basis of $\mR^m$.

Set
$$
\cE_0:=\sum^m_{i=1}q^2_i/\lambda_i,\ \ \ \cE_1:=\sum^m_{i=1}q^2_i.
$$
Then the quadratic variation of $\w(t)$ in $\mH^0$ and $\mH^1$ are given respectively  by
$$
[\w_\cdot]_{\mH^0}(t)=\cE_0 t,\ \ [\w_\cdot]_{\mH^1}(t)=\cE_1 t.
$$
We remark that $\cE_0\leq \cE_1/\l_1$.

Our main result in this section is the following:
\bt\label{Er}
Let $(\bT_t)_{t\geq 0}$ be the transition semigroup associated with (\ref{NSE00}).
For any sufficiently large $m_*=m_*(\cE_1,\l_1,N)\in\mN$,
there exists a unique invariant probability measure
associated with $(\bT_t)_{t\geq 0}$.
\et

We shall divide the proof into two parts. In the first part, we shall prove
the asymptotic strong Feller property of $(\bT_t)_{t\geq 0}$ (cf. \cite[Proposition 3.12]{Ha-Ma}).
In the second part, we shall prove a support property of the invariant measure, namely that
the origin $0$ is contained in the support of each invariant measure (cf. \cite{EM}).
By \cite[Proposition 3.12 and Corollary 3.17]{Ha-Ma}, these two parts will imply Theorem \ref{Er}.

\subsection{Asymptotic Strong Feller Property}
Let $\u(t,\om;\u_0)$ be the unique solution of Eq. (\ref{NSE00}).
For $0\leq s<t$,
let $\cJ_{s,t}$ denote the derivative flow of $\u(t,\om;\u_0)$ between $s$ and $t$ with respect to
the initial values $\u_0$,
i.e., for every $\v_0\in\mH^1$, $\cJ_{s,t}\v_0\in\mH^1$ satisfies
\begin{align}
\p_{t}\cJ_{s,t}\v_0=\Delta\cJ_{s,t}\v_0+\cK(\u(t,\om;\u_0),\cJ_{s,t}\v_0),\ \
\cJ_{s,s}\v_0=\v_0,\label{Der}
\end{align}
where $\cK$ is linear with respect to the second component and given by
$$
\cK(\u,\v):=-\sP((\v\cdot\nabla)\u+(\u\cdot\nabla)\v)-
\sP(g_N(|\u|^2)\v+2g'_N(|\u|^2)\<\u,\v\>_{\mR^3}\u).
$$
In the Appendix, we shall prove that for each $\om$
\begin{align}
(\cJ_{0,t}\v_0)(\om)=\lim_{\eps\downarrow 0}\frac{\u(t,\om;\u_0+\eps\v_0)-\u(t,\om;\u_0)}{\eps}
\quad \mbox{ in $\mH^1$}.\label{Der2}
\end{align}

Let us now consider the Malliavin derivative of $\u(t,\om;\u_0)$ with respect to $\om$.
Let $\sH$ be the Cameron-Martin space, i.e., all absolutely continuous functions from $\mR_+$
to $\mR^m$ with locally square integrable derivative. For any $v\in\sH$, the Malliavin derivative is
defined by
\begin{align}
D^v\u(t,\om;\u_0):=\lim_{\epsilon\rightarrow 0}
\frac{\u(t,\om+\epsilon v;\u_0)-\u(t,\om;\u_0)}{\epsilon}, \ \ P-a.s..\label{Lim1}
\end{align}
Notice that  $v$ can be random and possibly nonadapted to the filtration generated by $W$.
For the sake of simplicity, we write $\cA_t v:=D^v\u(t,\om;\u_0)$.
Then
\begin{align}
\p_{t}\cA_t v=\Delta\cA_t v+\cK(\u(t,\om;\u_0),\cA_t v)+Q\dot v(t),\ \ \cA_0 v=0,\label{Mal}
\end{align}
where $\dot v(t)$ is the derivative of $v(t)$ with respect to $t$.

By the formula of variation of constants, it is easy to see that
\begin{align*}
\cA_t v=\int^t_0\cJ_{s,t}Q\dot v(s)\dif s.
\end{align*}
Moreover, for any $\v_0\in\mH^1$ and $v\in\sH$, set
$$
\v(t):=\cJ_{0,t}\v_0-\cA_tv.
$$
Then
\begin{align}
\p_t \v(t)=\Delta\v(t)+\cK(\u(t),\v(t))-Q \dot v(t),\ \ \v(0)=\v_0.\label{IP1}
\end{align}

As done in \cite{Ha-Ma}, our main aim is to construct a suitable $v$ such that $\v(t)$
exponentially decays to zero in some sense as $t\rightarrow\infty$.
We first introduce some necessary notations and prove some preparing lemmas.

Let $\mH^1_\ell$ denote the following
finite dimensional subspace of $\mH^1$ (called low mode space)
$$
\mH^1_\ell:=\mbox{span}\{\e_1,\cdots, \e_m\}.
$$
Then we have the following direct sum decomposition:
$$
\mH^1=\mH^1_\ell\oplus\mH^1_\hbar
$$
and for any $\u\in\mH^1$,
$$
\u=\u_\ell+\u_\hbar, \ \ \ \u_\ell\in\mH^1_\ell, \ \u_\hbar\in\mH^1_\hbar.
$$
The co-dimensional space $\mH^1_\hbar$ is also called high mode space.
For any $\v\in\mH^0$, we define
$$
\Pi_\ell \v:=\sum_{i=1}^m\<(-\Delta)\e_i,\v\>_{\mH^0} \e_i\in\mH^1_\ell
$$
and
$$
\Pi_\hbar\v:=\v-\Pi_\ell\v\in\mH^0.
$$
In what follows, we shall always write $\v_\ell:=\Pi_\ell\v$ and $\v_\hbar:=\Pi_\hbar\v$.

The following lemma is immediate.
\bl\label{LE00}
For any $\u\in\mH^2$
$$
\|\Delta\Pi_\hbar\u\|^2_{\mH^0}\geq \lambda_m\|\nabla\Pi_\hbar\u\|^2_{\mH^0}.
$$
\el
We also need the following lemma. Recall that $\|\u\|_{\mH^2}=\|\Delta\u\|_{\mH^0}$.
\bl\label{LEE2}
For any $\u,\v\in\mH^2$, set
\begin{align}
\cN(\u):=\|\u\|^2_{\mH^2}+\||\u|\cdot|\nabla\u|\|_{L^2}^2.\label{cN}
\end{align}
Then
$$
\<\v_\hbar,\cK(\u,\v)\>_{\mH^1}\leq \frac{1}{2}\|\Delta\v_\hbar\|^2_{\mH^0}
+C_N\cdot\cN(\u)\cdot(\|\v_\hbar\|^2_{\mH^1}
+\|\v_\ell\|^2_{\mH^1})
$$
and
$$
\|\Pi_\ell\cK(\u,\v)\|^2_{\mH^1}\leq C_m\|\v\|^2_{\mH^0}\cdot(1+\|\u\|_{\mH^1}^4),
$$
where the constant $C_N$ (resp. $C_m$) only depends on $N$ (resp. $m$).
\el
\begin{proof}
For the first, we write
\begin{align*}
\<\v_\hbar,\cK(\u,\v)\>_{\mH^1}=I_1+I_2+I_3+I_4,
\end{align*}
where
\begin{align*}
I_1&:= -\<\v_\hbar,\sP((\v\cdot\nabla)\u)\>_{\mH^1},\\
I_2&:= -\<\v_\hbar,\sP((\u\cdot\nabla)\v)\>_{\mH^1},\\
I_3&:= -\<\v_\hbar,\sP(g_N(|\u|^2)\v)\>_{\mH^1},\\
I_4&:= -2\<\v_\hbar,\sP(g'_N(|\u|^2)\<\u,\v\>_{\mR^3}\u)\>_{\mH^1}.
\end{align*}
For $I_1$, by Young's inequality and the Sobolev inequality (\ref{Sob}) we have
\begin{align*}
I_1&\leq \frac{1}{8}\|\Delta\v_\hbar\|^2_{\mH^0}+2\||\v|\cdot|\nabla\u|\|^2_{L^2}
\leq \frac{1}{8}\|\Delta\v_\hbar\|^2_{\mH^0}+2\|\v\|^2_{L^6}\cdot\|\nabla\u\|^2_{L^3}\leq\\
&\leq \frac{1}{8}\|\Delta\v_\hbar\|^2_{\mH^0}+C\|\v_\hbar\|^2_{\mH^1}\cdot\|\u\|^2_{\mH^2}
+C\|\v_\ell\|^2_{\mH^1}\cdot\|\u\|^2_{\mH^2}.
\end{align*}
For $I_2$, we have
\begin{align*}
I_2&\leq \frac{1}{8}\|\Delta\v_\hbar\|^2_{\mH^0}+2\||\nabla\v|\cdot|\u|\|^2_{L^2}
\leq \frac{1}{8}\|\Delta\v_\hbar\|^2_{\mH^0}+2\|\v\|^2_{\mH^1}\cdot\|\u\|^2_{L^\infty}\leq\\
&\leq \frac{1}{8}\|\Delta\v_\hbar\|^2_{\mH^0}+C\|\v_\hbar\|^2_{\mH^1}\cdot\|\u\|^2_{\mH^2}
+C\|\v_\ell\|^2_{\mH^1}\cdot\|\u\|^2_{\mH^2}.
\end{align*}
For $I_3$, we have
\begin{align*}
I_3&= -\<\nabla\v_\hbar,g_N(|\u|^2)\nabla\v)\>_{\mH^0}
-\<\nabla\v_\hbar,g'_N(|\u|^2)\nabla|\u|^2\v)\>_{\mH^0}\\
&\leq \|\u\|^2_{L^\infty}\cdot\||\nabla\v_\hbar|\cdot|\nabla\v|\|_{L^1}
+\|\nabla\v_\hbar\|_{L^6}\cdot\|\v\|_{L^3}\cdot\|\nabla|\u|^2\|_{L^2}\\
&\leq C\|\u\|^2_{\mH^2}\cdot(\|\v_\hbar\|^2_{\mH^1}+\|\v_\ell\|^2_{\mH^1})
+\frac{1}{8}\|\v_\hbar\|^2_{\mH^2}
+C\|\v\|_{\mH^1}^2\cdot\||\u|\cdot|\nabla\u|\|_{L^2}^2.
\end{align*}
For $I_4$, noting that
$$
|g_N^{\prime\prime}(r)|\leq C \cdot 1_{\{N<r<N+1\}},
$$
we similarly have
$$
I_4\leq C\|\u\|^2_{\mH^2}\cdot(\|\v_\hbar\|^2_{\mH^1}+\|\v_\ell\|^2_{\mH^1})
+\frac{1}{8}\|\v_\hbar\|^2_{\mH^2}
+C_N\|\v\|_{\mH^1}^2\cdot\||\u|\cdot|\nabla\u|\|_{L^2}^2.
$$
Combining the above calculations, we obtain the first estimate.

As for the second one, we may write
$$
\|\Pi_\ell\cK(\u,\v)\|^2_{\mH^1}=\sum^m_{i=1}\<\e_i,\cK(\u,\v)\>_{\mH^1}^2
=\sum^m_{i=1}\left(\sum_{j=1}^4J_{ij}\right)^2,
$$
where
\begin{align*}
J_{i1}&:= -\<\e_i,\sP((\v\cdot\nabla)\u)\>_{\mH^1},\\
J_{i2}&:= -\<\e_i,\sP((\u\cdot\nabla)\v)\>_{\mH^1},\\
J_{i3}&:= -\<\e_i,\sP(g_N(|\u|^2)\v)\>_{\mH^1},\\
J_{i4}&:= -2\<\e_i,\sP(g'_N(|\u|^2)(\u\cdot\v)\u)\>_{\mH^1}.
\end{align*}
For $J_{i1}$, we have
\begin{align*}
J_{i1}= \<\Delta\e_i,(\v\cdot\nabla)\u)\>_{\mH^0}=-\<\nabla\Delta\e_i,\v\otimes\u\>_{\mH^0}
\leq \|\nabla\Delta\e_i\|_{L^\infty}\||\v|\cdot\u\|_{L^1}
\leq C_{\e_i}\|\v\|_{\mH_0}\|\u\|_{\mH_0}.
\end{align*}
Similarly, we have
\begin{align*}
J_{i2}&\leq C_{\e_i}\|\v\|_{\mH_0}\|\u\|_{\mH_0},\\
J_{i3}&\leq C_{\e_i}\|\v\|_{\mH_0}\|\u\|^2_{L^4},\\
J_{i4}&\leq C_{\e_i}\|\v\|_{\mH_0}\|\u\|^2_{L^4}.
\end{align*}
Summarizing the above calculations and by the Sobolev embedding theorem, we obtain the second estimate.
\end{proof}

We now prove the following crucial estimate about the solution $\u(t)$.
\bl\label{LEE3}
(i) For any $\eta>0$, there exist constants $C_\eta, C_{\cE_1,\l_1,N,\eta}>0$
such that for any $t>0$ and $\u_0\in\mH^1$
$$
\mE\exp\left\{\eta\int^t_0\cN(\u(s;\u_0))\dif s\right\}
\leq \exp\{C_\eta \|\u_0\|^2_{\mH^1}+C_{\cE_1,\l_1,N,\eta}t\},
$$
where $\cN(\u)$ is defined by (\ref{cN}).

(ii) There exist constants $C_N, C_{\cE_1,\l_1,N}>0$ such that for any $t>0$ and $\u_0\in\mH_1$
$$
\mE\|\u(t;\u_0)\|^2_{\mH^1}\leq \|\u_0\|^2_{\mH^1}(C_N\cdot t+1)e^{-t/2}+
C_{\cE_1,\l_1,N}.
$$
\el
\begin{proof}
By It\^o's formula, we have
\begin{align}
\dif \|\u(t)\|_{\mH^0}^2=2\<\u(t), A(\u(t))\>_{\mH^0}\dif t
+2\<\u(t),\dif \w(t)\>_{\mH^0}+\cE_0\dif t.   \label{PO1}
\end{align}
By (\ref{Es44}) and Young's inequality, we know
\begin{align}
\<\u(t), A(\u(t))\>_{\mH^0}&\leq -\|\nabla\u(t)\|^2_{\mH^0}-\|\u(t)\|^4_{L^4}+N\|\u(t)\|^2_{\mH^0}\no\\
&\leq -\|\u(t)\|^2_{\mH^1}-\frac{1}{2}\|\u(t)\|^4_{\mH^0}+\frac{N^2}{2}.\label{PO2}
\end{align}
Using Lemma \ref{LEE1} in the Appendix, we get for any $t,\eta>0$
\begin{align}
\mE\exp\left\{\eta\int^t_0\|\u(s)\|_{\mH^1}^2\dif s\right\}\leq \exp\{\eta \|\u_0\|^2_{\mH^0}
+C_{\cE_0,N,\eta}t\}.\label{Ep2}
\end{align}
Again, by It\^o's formula  and (\ref{Es4}), we have
\begin{align}
\dif \|\u(t)\|_{\mH^1}^2&= 2\lb\u(t), A(\u(t))\rb\dif t+2\<\u(t),\dif \w(t)\>_{\mH^1}
+\cE_1\dif t\no\\
&\leq (-\cN(\u(t))+C_N\|\u(t)\|^2_{\mH^1})\dif t+2\<\u(t),\dif \w(t)\>_{\mH^1}+\cE_1\dif t.\label{PPO1}
\end{align}
As in the proof of Lemma \ref{LEE1} in the Appendix, using (\ref{Ep2}) and exponential martingales,
we then get the first estimate.

On the other hand, from (\ref{PO1}) and (\ref{PO2}), we have
\begin{align*}
\dif \|\u(t)\|_{\mH^0}^2\leq-\frac{1}{2}\|\u(t)\|^2_{\mH^0}\dif t
+2\<\u(t),\dif \w(t)\>_{\mH^0}+(\cE_0+N^2+\frac{1}{2})\dif t.
\end{align*}
It is direct by Gronwall's inequality that
\begin{align*}
\mE\|\u(t;\u_0)\|^2_{\mH^0}\leq \|\u\|^2_{\mH^0} e^{-t/2}+2(\cE_0+N^2+1).
\end{align*}
Thus, thanks to
\begin{align}
\|\u(t)\|^2_{\mH^1}\leq \|\u(t)\|_{\mH^2}\|\u(t)\|_{\mH^0},\label{PP5}
\end{align}
we obtain
\begin{align*}
\dif (e^{t/2}\|\u(t)\|_{\mH^1}^2)&= e^{t/2}(2\lb\u(t), A(\u(t))\rb+\cE_1
+\frac{1}{2}\|\u(t)\|_{\mH^1}^2)\dif t+2e^{t/2}\<\u(t),\dif \w(t)\>_{\mH^1}\\
&\leq e^{t/2}(-\|\u(t)\|^2_{\mH^2}+C_N\|\u(t)\|^2_{\mH^1}+\cE_1)\dif t
+2e^{t/2}\<\u(t),\dif \w(t)\>_{\mH^1}\\
&\leq e^{t/2}(C_{N}\|\u(t)\|^2_{\mH^0}+\cE_1)\dif t
+2e^{t/2}\<\u(t),\dif \w(t)\>_{\mH^1}.
\end{align*}
Therefore,
\begin{align*}
 e^{t/2}\mE\|\u(t;\u_0)\|_{\mH^1}^2&\leq \|\u_0\|^2_{\mH^1}
+C_N\int^t_0e^{s/2}\mE\|\u(s;\u_0)\|^2_{\mH^0}\dif s+2\cE_1 e^{t/2}\\
&\leq \|\u_0\|^2_{\mH^1}+C_N\int^t_0(\|\u_0\|^2_{\mH^0}
+2e^{s/2}(\cE_0+N^2+1))\dif s+2\cE_1 e^{t/2}\\
&\leq \|\u_0\|^2_{\mH^1}+C_N \|\u_0\|^2_{\mH^0}t+
e^{t/2}(C_N(\cE_0+N^2+1)+2\cE_1),
\end{align*}
which then gives the second estimate.
\end{proof}

Based on the previous discussions and lemmas, we can now prove the following proposition,
which will imply the asymptotic strong Feller property
of $(\bT_t)_{t\geq 0}$ according to \cite[Proposition 3.12]{Ha-Ma}.
\bp\label{Pro1}
Let $(\bT_t)_{t\geq 0}$ be the semigroup associated with (\ref{NSE00}).
There exist a constant $m_*:=m_*(\cE_1,N)\in\mN$ and
constants $C_0,C_1,\gamma>0$ such that for any $t>0$, $\u_0\in\mH^1$,
and any Fr\'echet differentiable function $\varphi$ on $\mH^1$ with
$\|\varphi\|_\infty,\|\nabla\varphi\|_\infty<+\infty$,
$$
\|\nabla\bT_t\varphi(\u_0)\|_{\mH^1}\leq C_0\cdot\exp\{C_1\|\u_0\|^2_{\mH^1}\}
\cdot(\|\varphi\|_\infty+e^{-\gamma t}\|\nabla\varphi\|_\infty).
$$
\ep
\begin{proof}
For any $\v_0\in\mH^1$ with $\|\v_0\|_{\mH^1}=1$, define
\begin{align*}
\v_\ell(t):=\left\{
\begin{array}{ll}
\v_{0\ell}\cdot(1-t/(2\|\v_{0\ell}\|_{\mH^1})),& t\in[0,2\|\v_{0\ell}\|_{\mH^1}]\\
0,&t\in( 2\|\v_{0\ell}\|_{\mH^1},\infty).
\end{array}
\right.
\end{align*}
Let $\v_\hbar(t)$ solve the following linear evolution equation:
\begin{align*}
\p_t \v_\hbar(t)=\Delta\Pi_\hbar\v_\hbar(t)+\Pi_\hbar\cK(\u(t),\v_\hbar(t)+\v_\ell(t)),\ \
\v_\hbar(0)=\v_{0\hbar}.
\end{align*}

Set
$$
\v(t):=\v_\ell(t)+\v_\hbar(t)
$$
and
$$
\dot v(t):=Q^{-1}\left(\frac{\v_\ell\cdot 1_{\{t<2\|\v_\ell\|_{\mH^1}\}}}{2\|\v_\ell\|_{\mH^1}}
+\Delta\v_\ell(t)+\Pi_\ell\cK(\u(t),\v(t))\right).
$$
Then $v(t)\in\sH$ is a continuous adapted process.
From the construction, one finds that $\v(t)$ together with $v(t)$ solves the equation (\ref{IP1}).

Thus, we have
\begin{align}
\<\nabla\bT_t\varphi(\u_0),\v_0\>_{\mH^1}&=
\mE\<(\nabla\varphi)(\u(t;\u_0)),\cJ_{0,t}\v_0\>_{\mH^1}\no\\
&= \mE\<(\nabla\varphi)(\u(t;\u_0)),\cA_{t}v(t)\>_{\mH^1}+
\mE\<(\nabla\varphi)(\u(t;\u_0)),\v(t)\>_{\mH^1}\no\\
&= \mE(D^v(\varphi(\u(t;\u_0))))+
\mE\<(\nabla\varphi)(\u(t;\u_0)),\v(t)\>_{\mH^1}\no\\
&= \mE\left(\varphi(\u(t;\u_0))\cdot\int^t_0\dot v(s)\dif W_s\right)+
\mE\<(\nabla\varphi)(\u(t;\u_0)),\v(t)\>_{\mH^1}\no\\
&\leq \|\varphi\|_\infty\left(\int^t_0\mE|\dot v(s)|^2\dif s\right)^{1/2}+
\|\nabla\varphi\|_\infty\mE\|\v(t)\|_{\mH^1},\label{LP4}
\end{align}
where the last equality is due to the integration by parts
formula in the Malliavin calculus (cf. \cite{M}).

By the chain rule and Lemmas \ref{LEE2} and \ref{LE00}, we have
\begin{align*}
\p_t\|\v_\hbar(t)\|^2_{\mH^1}&= -2\|\Delta\Pi_\hbar\v_\hbar(t)\|^2_{\mH^0}+
2\<\v_\hbar(t),\Pi_\hbar\cK(\u(t),\v(t))\>_{\mH^1}\\
&\leq -\|\Delta\Pi_\hbar\v_\hbar(t)\|^2_{\mH^0}+
C_N\cdot\cN(\u(t))\cdot(\|\v_\hbar(t)\|^2_{\mH^1}+\|\v_\ell(t)\|^2_{\mH^1})\\
&\leq (-\l_m+C_N\cdot\cN(\u(t)))\cdot\|\v_\hbar(t)\|^2_{\mH^1}+
C_N\cdot\cN(\u(t))\cdot\|\v_\ell(t)\|^2_{\mH^1}.
\end{align*}
Noting that $v(t)=0$ for $t\geq 2$, by Gronwall's inequality we get
\begin{align*}
\|\v_\hbar(t)\|^2_{\mH^1}&\leq \|\v_\hbar(0)\|^2_{\mH^1}
\exp\left\{-\l_m t+C_N\int^t_0\cN(\u(s))\dif s\right\}\\
&\quad +\exp\left\{-\l_m (t-2)+C_N\int^t_0\cN(\u(s))\dif s\right\}
\int^2_0\|\v_\ell(s)\|^2_{\mH^1}\dif s.
\end{align*}
By (i) of Lemma \ref{LEE3}, since $\l_m\uparrow\infty$ as $m\rightarrow\infty$,
there exist constants $\gamma>0$ and $m_*=m_*(\cE_1,\l_1,N)\in\mN$ such that for all $t\geq 0$,
$$
\mE\|\v_\hbar(t)\|^4_{\mH^1}\leq C_{\cE_1,\l_1,N}\cdot e^{C_N\|\u_0\|^2_{\mH^1}-\gamma t}.
$$
Hence, for any $t\geq 2$,
\begin{align}
\mE\|\v(t)\|_{\mH^1}\leq C_{\cE_1,\l_1,N}\cdot e^{C_N\|\u_0\|^2_{\mH^1}-\gamma t}.\label{LP3}
\end{align}
On the other hand, by  Lemma \ref{LEE2}, we have
\begin{align}
\mE|\dot v(t)|^2&\leq C_m\left(1+\mE(\|\v(t)\|^2_{\mH^0}(1+\|\u(t)\|_{\mH^1}^4))\right)\no\\
&\leq C_m\left(1+(\mE\|\v(t)\|^4_{\mH^0})^{1/2}(1+\mE\|\u(t)\|_{\mH^1}^8)^{1/2}\right).\label{PL1}
\end{align}
Using It\^o's formula and (\ref{Es444}), as in the proof of Theorem \ref{Th3}, we have
$$
\mE\|\u(t)\|^{2p}_{\mH^0}\leq C\|\u_0\|^{2p}_{\mH^0}(1+t),
$$
and also by (\ref{Es4}), $\|\u\|_{\mH^1}^2\leq\|\u\|_{\mH^0}\|\u\|_{\mH^2}$ and Young's inequality,
\begin{align*}
\mE\|\u(t)\|^{2p}_{\mH^1}&\leq \|\u_0\|^{2p}_{\mH^1}-p\mE\int^t_0\|\u(s)\|^{2(p-1)}_{\mH^1}\|\u(s)\|^2_{\mH^2}\dif s
+C_N\mE\int^t_0\|\u(s)\|^{2p}_{\mH^1}\dif s+Ct\\
&\leq \|\u_0\|^{2p}_{\mH^1}-\frac{p}{2}\mE\int^t_0\|\u(s)\|^{2(p-1)}_{\mH^1}\|\u(s)\|^2_{\mH^2}\dif s
+C_N\mE\int^t_0 \|\u(s)\|^{2(p-1)}_{\mH^1}\|\u(s)\|^2_{\mH^0}\dif s+Ct\\
&\leq \|\u_0\|^{2p}_{\mH^1}+C_N\mE\int^t_0 \|\u(s)\|^{2p}_{\mH^0}\dif s+Ct\leq C\|\u_0\|^{2p}_{\mH^1}(1+t^2).
\end{align*}
Thus, integrating both sides of (\ref{PL1}) and using (\ref{LP3}), we obtain
\begin{align}
\int^\infty_0\mE|\dot v(t)|^2\dif t\leq C_{m,\cE_1,\l_1,N,\gamma}\cdot e^{C_N\|\u_0\|^2_{\mH^1}}\cdot
\left(1+\int^\infty_0 e^{-\gamma t}(1+t)\dif  t\right)
\leq C_{m,\cE_1,\l_1,N,\gamma}\cdot e^{C_N\|\u_0\|^2_{\mH^1}}.\label{LP5}
\end{align}
The proof is thus completed by combining (\ref{LP4})-(\ref{LP5}).
\end{proof}
\subsection{A Support Property of Invariant Measures}
\bp\label{Lp6}
The point $0$ belongs to the support of any invariant measure of $(\bT_t)_{t\geq 0}$.
\ep
For the proof we need the following lemma, whose proof in turn
is inspired by \cite{EM}.
\bl\label{LEE4}
For any $r_1,r_2>0$, there exists $T>0$ such that
\begin{align*}
\inf_{\|\u_0\|_{\mH^1}\leq r_1}P\{\om: \|\u(T,\om;\u_0)\|_{\mH^1}\leq r_2\}>0.
\end{align*}
\el
\begin{proof}
Set
$$
\v(t):=\u(t)-\w(t).
$$
Then
$$
\v'(t)=A(\v(t)+\w(t)),\ \ \v(0)=\u_0.
$$

Let $T>0$ and $\epsilon\in(0,1)$, to be determined below. We assume that
\begin{align}
\sup_{t\in[0,T]}\|\w(t)\|_{\mH^6}<\epsilon.\label{PP4}
\end{align}

First of all, by  the chain rule, we have
$$
\frac{\dif}{\dif t}\|\v(t)\|^2_{\mH^0}=J_1+J_2+J_3+J_4,
$$
where
\begin{align*}
J_1&:= -2\|\nabla\v(t)\|^2_{\mH^0}+2\<\Delta\w(t),\v(t)\>_{\mH^0},\\
J_2&:= -2\<\v(t),((\v(t)+\w(t))\cdot\nabla)(\v(t)+\w(t))\>_{\mH^0},\\
J_3&:= -2\<\v(t)+\w(t),g_N(|\v(t)+\w(t)|^2)(\v(t)+\w(t))\>_{\mH^0},\\
J_4&:= 2\<\w(t),g_N(|\v(t)+\w(t)|^2)(\v(t)+\w(t))\>_{\mH^0}.
\end{align*}

For $J_1$, by (\ref{PP4}) we have
\begin{align*}
J_1\leq -2\|\nabla\v(t)\|^2_{\mH^0}+C\epsilon\|\v(t)\|_{\mH^0}.
\end{align*}
Here and below, $C$ denotes an absolute constant.

For $J_2$, by the Sobolev inequality (\ref{Sob}) and (\ref{PP4}) we have
\begin{align*}
J_2&= -2\<\w(t),((\v(t)+\w(t))\cdot\nabla)(\v(t)+\w(t))\>_{\mH^0}\\
&\leq 2\|\nabla\w(t)\|_{L^\infty}\|\v(t)+\w(t)\|^2_{\mH^0}\\
&\leq C\epsilon\cdot\|\v(t)\|^2_{\mH^0}+C\epsilon.
\end{align*}

For $J_3$, we obviously have
$$
J_3\leq 0.
$$

For $J_4$, by (\ref{Sob}) and Young's inequality we have
\begin{align*}
J_4&\leq 2\|\w(t)\|_{L^\infty}\cdot\|\v(t)+\w(t)\|^3_{L^3}\leq C\epsilon\cdot\|\v(t)\|^3_{L^3}+C\epsilon^4\leq\\
&\leq C\epsilon\cdot\|\nabla\v(t)\|^{3/2}_{\mH^0}\|\v(t)\|^{3/2}_{\mH^0}+C\epsilon^4\\
&\leq \|\nabla\v(t)\|^2_{\mH^0}+
C\epsilon^4\cdot\|\v(t)\|^6_{\mH^0}+C\epsilon^4.
\end{align*}

Combing the above calculations gives that
\begin{align*}
\frac{\dif}{\dif t}\|\v(t)\|^2_{\mH^0}&\leq -\|\nabla\v(t)\|^2_{\mH^0}
+C\epsilon\cdot\|\v(t)\|^6_{\mH^0}+C\epsilon\\
&\leq -\frac{1}{\l_1}\|\v(t)\|^2_{\mH^0}
+C\epsilon\cdot\|\v(t)\|^6_{\mH^0}+C\epsilon,
\end{align*}
where the second step is due to the Poincare inequality (\ref{Poi}).

Note that $\|\v(t)\|^2_{\mH^0}$ depends on $\eps$ through (\ref{PP4}).
By Lemma \ref{Le10} in the Appendix, for any $\delta,h>0$,
we may choose a $T_0>0$ sufficiently large and
an $\epsilon$ small enough such that
\begin{align}
\sup_{t\in[0,T_0]}\|\v(t)\|_{\mH^0}\leq 2r_1\label{Po1}
\end{align}
and
\begin{align}
\sup_{t\in[T_0,T_0+h]}\|\v(t)\|_{\mH^0}<\delta.\label{Po2}
\end{align}

Let us now turn to the  estimate of the first order Sobolev norm of $\v(t)$.
By the chain rule again, we have
\begin{align*}
\frac{\dif}{\dif t}\|\v(t)\|^2_{\mH^1}=2\lb\v(t),A(\v(t)+\w(t))\rb=I_1+I_2+I_3+I_4,
\end{align*}
where
\begin{align*}
I_1&:= 2\lb\v(t)+\w(t),A(\v(t)+\w(t))\rb,\\
I_2&:= -2\<\Delta^2\w(t),\v(t)+\w(t)\>_{\mH^0},\\
I_3&:= 2\<\Delta\w(t),((\v(t)+\w(t))\cdot\nabla)(\v(t)+\w(t))\>_{\mH^0},\\
I_4&:= 2\<\Delta\w(t),g_N(|\v(t)+\w(t)|^2)(\v(t)+\w(t))\>_{\mH^0}.
\end{align*}

For $I_1$, by (\ref{Es4}) and (\ref{PP5}) we have
\begin{align*}
I_1&\leq -\|\v(t)+\w(t)\|^2_{\mH^2}+C_N\|\v(t)+\w(t)\|^2_{\mH^1}\\
&\leq -\frac{1}{2}\|\v(t)\|^2_{\mH^2}+4N\|\v(t)\|^2_{\mH^1}+C_N\epsilon\\
&\leq -\frac{1}{4}\|\v(t)\|^2_{\mH^2}+C_N\|\v(t)\|^2_{\mH^0}+C_N\epsilon.
\end{align*}
Here and below, $C_N$ denotes a constant only depending on $N$.

For $I_2$, we have
$$
I_2\leq C\epsilon+C\epsilon\|\v(t)\|_{\mH^0}.
$$

For $I_3$, we have
\begin{align*}
I_3\leq 2\|\nabla\Delta\w(t)\|_{L^\infty}\|\v(t)+\w(t)\|^2_{\mH^0}\leq C\epsilon+C\epsilon\|\v(t)\|^2_{\mH^0}.
\end{align*}

For $I_4$, we have
$$
I_4\leq C\epsilon+C\epsilon\|\v(t)\|^3_{L^3}\leq \frac{1}{8}\|\v(t)\|^2_{\mH^2}
+C\epsilon+C\epsilon\|\v(t)\|^6_{\mH^0}.
$$

Combing the above calculations gives that
\begin{align*}
\frac{\dif}{\dif t}\|\v(t)\|^2_{\mH^1}&\leq -\frac{1}{8}\|\v(t)\|^2_{\mH^2}
+C_N\cdot\|\v(t)\|^2_{\mH^0}+C\epsilon\cdot\|\v(t)\|^6_{\mH^0}+C\epsilon\\
&\leq -C_0\|\v(t)\|^2_{\mH^1}+C_N\cdot\|\v(t)\|^2_{\mH^0}+C\epsilon\|\v(t)\|^6_{\mH^0}+C\epsilon.
\end{align*}

By Gronwall's inequality, for any $0<t_1<t_2$ we have
\begin{align*}
\|\v(t_2)\|^2_{\mH^1}&\leq e^{-C_0(t_2-t_1)}\|\v(t_1)\|^2_{\mH^1}
+\frac{1}{C_0}\Big(C_N\cdot\sup_{t\in[t_1,t_2]}\|\v(t)\|^2_{\mH^0}
+C\epsilon\cdot\sup_{t\in[t_1,t_2]}\|\v(t)\|^6_{\mH^0}+C\epsilon\Big).
\end{align*}

Firstly, letting $t_1=0$ and $t_2=T_0$  and by (\ref{Po1}), we find
\begin{align*}
\|\v(T_0)\|^2_{\mH^1}\leq r_1^2
+\frac{1}{C_0}\Big(C_N\cdot\sup_{t\in[0,T_0]}\|\v(t)\|^2_{\mH^0}
+C\epsilon\cdot\sup_{t\in[0,T_0]}\|\v(t)\|^6_{\mH^0}+C\epsilon\Big)\leq C_{N,C_0}(r_1^6+1).
\end{align*}

Secondly, letting $t_1=T_0$ and $t_2=T_0+h$ yields
\begin{align*}
\|\v(T_0+h)\|^2_{\mH^1}\leq e^{-C_0h}C_{N,C_0}(r_1^6+1)
+\frac{1}{C_0}\Big(C_N\cdot\sup_{t\in[T_0,T_0+h]}\|\v(t)\|^2_{\mH^0}
+C\epsilon\cdot\sup_{t\in[T_0,T_0+h]}\|\v(t)\|^6_{\mH^0}+C\epsilon\Big),
\end{align*}
which together with (\ref{Po2}) implies that for some $T$ large enough and $\epsilon>0$ small enough
$$
\|\v(T)\|_{\mH^1}\leq r_2/2.
$$
Therefore, there exist  $T$ sufficiently large and $\epsilon$ small enough such that
for any $\|\u_0\|_{\mH^1}\leq r_1$
$$
\|\u(T,\om;\u_0)\|_{\mH^1}\leq r_2.
$$
That is, if we set
$$
\Omega_\epsilon:=\left\{\om: \sup_{t\in[0,T]}\|\w(t,\om)\|_{\mH^6}<\epsilon\right\},
$$
then
$$
\Omega_\epsilon\subset\cap_{\|\u_0\|_{\mH^1}\leq r_1}\{\om: \|\u(T,\om;\u_0)\|_{\mH^1}\leq r_2\}.
$$
The desired estimate now follows from the fact that $\Omega_\epsilon$ is an open subset
of $\Omega$ and
$P(\Omega_\epsilon)>0$.
\end{proof}

\vspace{4mm}

{\it Proof of Proposition \ref{Lp6}}:
For $r>0$, let $\mB_r:=\{\u_0\in\mH^1: \|\u_0\|_{\mH^1}\leq r\}$ be the ball in $\mH^1$.
For each invariant measure $\mu$, we can choose some $r_1>0$ such that
$$
\mu(\mB_{r_1})\geq 1/2.
$$
By Lemma \ref{LEE4}, we further have for any $r_2>0$ and some $t>0$
\begin{align*}
\mu(\mB_{r_2})= \int_{\mH^1}(\bT_t1_{\mB_{r_2}})(\u_0)\mu(\dif\u_0)
\geq\int_{\mB_{r_1}}(\bT_t1_{\mB_{r_2}})(\u_0)\mu(\dif\u_0)
\geq \mu(\mB_{r_1})\cdot \inf_{\u_0\in \mB_{r_1}}(\bT_t1_{\mB_{r_2}})(\u_0) >0,
\end{align*}
which means that $0$ belongs to the support of $\mu$.

\vspace{4mm}

{\it Proof of Theorem \ref{Er}}: The assertion follows from
Propositions \ref{Pro1} and \ref{Lp6} due to
\cite[Proposition 3.12, Corollary 3.17]{Ha-Ma}.

\section{Appendix}

\subsection{Proof of Proposition \ref{Th2}}
In this subsection, we prove the martingale characterization of weak solutions.

First of all, (i)$\Longrightarrow$(ii) is direct by It\^o's formula.
Let us prove (ii)$\Longrightarrow$(i). Define for $\e\in\sE$ (see Subsection 2.3
for the notation $\sE$)
$$
M_{\e}(t,\u):=\<\u(t)-\u(0),\e\>_{\mH^1}-\int^t_0\lb A(\u(s)),\e\rb\dif s
-\int^t_0\<\f(s,\u(s)),\e\>_{\mH^1}\dif s.
$$
Using (ii) and by  simple approximations as in \cite{St}, one knows that $\{M_{\e}(t,\u), t\geq 0\}$
is a continuous local martingale under $P_\vartheta$ with respect to $\cB_t(\mX)$, and its quadratic
variation process is given by
$$
[M_{\e}](t,\u)=\int^t_0\|\<B(s,\u(s)),\e\>_{\mH^1}\|^2_{l^2}\dif s.
$$

Set
\begin{align}
M(t,\u):=\sum_{j=1}^\infty M_{\e_j}(t,\u)\e_j.\label{Se}
\end{align}
Then $t\mapsto M(t,\u)$ is an $\mH^1$-valued continuous local
martingale under $P_\vartheta$ with respect to $\cB_t(\mX)$. Indeed, for any $R>0$,
define the stopping time
$$
\tau_R(\u):=\inf\left\{t\geq 0: \int^t_0\|B(s,\u(s))\|^2_{L_2(l^2;\mH^1)}\dif s\geq R\right\}.
$$
Then by (\ref{Es00}) and (\ref{LP2}) we have
$$
\tau_R(\u)\uparrow\infty,\ \ P_\vartheta-a.a.\ \u,  \ \mbox{ as $R\to\infty$.}
$$
Set
$$
M^{R,n}(t,\u):=\sum_{j=1}^n M_{\e_j}(t\wedge\tau_R,\u)\e_j.
$$
It is clear that $M^{R,n}(t,\u)$ is an $\mH^1$-valued continuous martingale with
\begin{align*}
\ll M^{R,n}\gg_{\mH^1}(t,\u)&=\sum_{i,j=1}^n[M_{\e_i},M_{\e_j}](t\wedge\tau_R,\u)\cdot \e_i\otimes \e_j\\
&= \sum_{i,j=1}^n\int^{t\wedge\tau_R}_0\<\<B(s,\u(s)),\e_j\>_{\mH^1},
\<B(s,\u(s)),\e_j\>_{\mH^1}\>_{l^2}\cdot  \e_i\otimes \e_j \dif s,
\end{align*}
where $\ll \cdot\gg_{\mH^1}$ denotes the square variation of $M^R$ in $\mH^1$.
Moreover, by Burkholder's inequality we have, for any $T>0$
\begin{align*}
\mE^{P_\vartheta}\left(\sup_{t\in[0,T]}\|M^{R,n}(t,\u)-M^{R,m}(t,\u)\|_{\mH^1}^2\right)
\leq C\sum_{j=n}^m\mE^{P_\vartheta}\left(\int^{T\wedge\tau_R}_0
\|\<B(s,\u(s)),\e_j\>_{\mH^1}\|_{l^2}^2\dif s\right)\rightarrow  0
\end{align*}
as $n,m\rightarrow\infty$. Hence, the series in (\ref{Se}) converges in $C([0,T];\mH^1)$, $P_\vartheta$-a.s.,
and $M^R(t,\u):=M(t\wedge\tau_R,\u)$ is
an $\mH^1$-valued  continuous  square integrable martingale with
\begin{align*}
\ll M^R\gg_{\mH^1}(t,\u)&=\sum_{i,j=1}^\infty[M_{\e_i},M_{\e_j}](t\wedge\tau_R,\u)\cdot \e_i\otimes \e_j\\
&= \sum_{i,j=1}^\infty\int^{t\wedge\tau_R}_0\<\<B(s,\u(s)),\e_j\>_{\mH^1},
\<B(s,\u(s)),\e_j\>_{\mH^1}\>_{l^2}\cdot  \e_i\otimes \e_j \dif s.
\end{align*}
Letting $R\to \infty$ we obtain the desired property of $M(t,\u)$.

In particular, the following equality holds in $\mH^0$
$$
\u(t)=\u(0)+\int^t_0A(\u(s))\dif s+\int^t_0\f(s,\u(s))\dif s+M(t,\u),\ \ P_\vartheta-a.s..
$$
By It\^o's formula  (cf. \cite{Ro, Roe}), we obtain that $P_\vartheta(C([0,\infty),\mH^1))=1$.
The existence of weak solutions now follows from the
representation theorem for martingales (cf. \cite[Lemma 3.2]{Mi-Ro1} or \cite[Theorem 8.2]{DaZa}).

\subsection{Two Basic Estimates}
In this subsection, we prove two basic estimates used in Section 5.
\bl\label{Le10}
Let $\{\varphi_\epsilon(\cdot,r_0),\epsilon\in(0,1),r_0\geq0\}$ be a family of positive
real functions on $\mR_+$ with $\varphi_\epsilon(0,r_0)=r_0$.
Suppose that for some $p>1$, $C_0,C_1,C_2>0$, $C_3\geq 0$ and any $\epsilon\in(0,1)$ and $t\geq 0$
$$
\varphi'_\epsilon(t,r_0)\leq -C_0\varphi_\epsilon(t,r_0)
+C_1\epsilon\cdot\varphi_\epsilon(t,r_0)^p+C_2\epsilon+C_3.
$$
Then: (i) For any $T>0$ and $R>0$, there exists  $\epsilon_0>0$ such that
$$
\sup_{t\in[0,T],\epsilon\in[0,\epsilon_0],r_0\in[0,R]}\varphi_\epsilon(t,r_0)\leq 2R+2C_3/C_0.
$$

(ii) If $C_3=0$, then for any $\delta>0$ and $R,h>0$, there exist  $T>0$ and $\epsilon_0>0$ such that
$$
\sup_{t\in[T,T+h],\epsilon\in[0,\epsilon_0],r_0\in[0,R]}\varphi_\epsilon(t,r_0)\leq \delta.
$$
\el
\begin{proof}
Let $C^\epsilon_4:=(C_2\epsilon+C_3)/C_0$ and set
$$
\phi(t):=e^{C_0 t}(\varphi_\epsilon(t,r_0)-C_4^\epsilon).
$$
Then for fixed $T>0$ and any $t\in[0,T]$
\begin{align*}
\phi'(t)\leq C_1 e^{C_0 t}\epsilon\cdot\varphi_\epsilon(t,r_0)^p
\leq C_1 \epsilon\cdot(\phi(t)+C_4^\epsilon \cdot e^{C_0T})^p.
\end{align*}
Solving this differential inequality gives that
\begin{align*}
\phi(T)\leq \Big[(\phi(0)+C_4^\epsilon \cdot e^{C_0T})^{1-p}
+C_1(1-p)\epsilon T\Big]^{\frac{1}{1-p}}-C_4^\epsilon \cdot e^{C_0T}.
\end{align*}
Hence,
\begin{align*}
\varphi_\epsilon(T,r_0)&\leq e^{-C_0T}\Big[(r_0+C_4^\epsilon \cdot (e^{C_0T}-1))^{1-p}+C_1(1-p)\epsilon T\Big]^{\frac{1}{1-p}}\\
&\leq \Big[(e^{-C_0T}R+C_4^\epsilon)^{1-p}+C_1(1-p)\epsilon Te^{(p-1)C_0T}\Big]^{\frac{1}{1-p}}.
\end{align*}
Now the assertions easily follow by suitable choices of $\epsilon$ and $T$.
\end{proof}

We now prove the following exponential estimate.
\bl\label{LEE1}
Let $X_t$ be a positive It\^o process of the form
\begin{align}
X_t=x_0+\int^t_0M_s\dif W_s+\int^t_0N_s\dif s,\label{ML4}
\end{align}
where $s\mapsto M_s, N_s$ are two measurable adapted processes.
Suppose that there exist a positive process $Y_s$ and some $\alpha>1$ and $C_0,C_1,C_2,C_3>0$
such that for any $s\geq 0$
\begin{align}
N_s\leq -C_0 X^\a_s-Y_s+C_1,\ \ \ |M_s|^2\leq C_2X_s+C_3.\label{ML1}
\end{align}
Then for any $t,\eta>0$
\begin{align}
\mE e^{\eta X_t}\leq C_{\a,\eta}\cdot \exp\{ e^{-C_0t/2} \eta x_0\}\label{ML2}
\end{align}
and
\begin{align}
\mE \exp\left\{\eta \int^t_0(\frac{C_0}{2}X^\a_s+Y_s)\dif s\right\}
\leq\exp\{\eta x_0+C_{\a,\eta} t\}.\label{ML3}
\end{align}
\el
\begin{proof}
Let us first prove that for any $t,\eta>0$
\begin{align}
\mE e^{\eta X_t}<+\infty.\label{Ep1}
\end{align}
Set for $R>0$
$$
\tau_R:=\inf\{t\geq 0: |X_t|\geq R\}.
$$
By It\^o's formula, (\ref{ML1}) and Young's inequality, we have
\begin{align*}
\dif e^{\eta X_t}&= \eta e^{\eta X_t}M_t\dif W_t+
\eta e^{\eta X_t}N_t\dif t+\frac{\eta^2}{2}e^{\eta X_t}|M_t|^2\dif t\\
&\leq \eta e^{\eta X_t}M_t\dif W_t+\eta e^{\eta X_t}(-C_0 X_t^\a+C_1
+\frac{\eta}{2}(C_2X_t+C_3))\dif t\\
&\leq \eta e^{\eta X_t}M_t\dif W_t+\eta e^{\eta X_t}(-\frac{C_0}{2} X_t^\a
+C_{\a,\eta})\dif t\\
&\leq \eta e^{\eta X_t}M_t\dif W_t+\eta e^{\eta X_t}(-\frac{C_0}{2} X_t
+C_{\a,\eta})\dif t.
\end{align*}
Set
$$
f_R(t):=\mE e^{\eta X_{t\wedge\tau_R}}.
$$
Then
$$
f'_R(t)\leq C_{\a,\eta}f_R(t).
$$
Hence
\begin{align*}
f_R(t)=\mE e^{\eta X_{t\wedge\tau_R}}\leq e^{\eta x_0} e^{C_{\a,\eta} t}.
\end{align*}
By Fatou's lemma, we obtain (\ref{Ep1}).

We now set
$$
f(t):=\mE e^{\eta X_{t}}.
$$
Then by Jensen's inequality, we obtain
$$
f'(t)\leq-\frac{C_0}{2}f(t)\log f(t)+C_{\a,\eta} f(t).
$$
Solving this differential equality gives the first estimate (\ref{ML2}).

On the other hand, for any $t,\eta>0$, we have by (\ref{ML4}) and (\ref{ML1})
\begin{align*}
\eta\int^t_0(\frac{C_0}{2}X^\a_s+Y_s)\dif s
\leq \eta x_0+\eta\int^t_0M_s\dif W_s+\int^t_0(-\frac{\eta C_0}{2}X^\a_s+C_1\eta)\dif s.
\end{align*}
Noting that by (\ref{ML1}) and (\ref{ML2})
\begin{align*}
\mE \exp\left\{\frac{\eta^2}{2}\int^t_0|M_s|^2\dif s\right\}
\leq\frac{1}{t}\int^t_0\mE\exp\{t\eta^2|M_s|^2/2\}\dif s<+\infty,
\end{align*}
we know by Novikov's criterion that
\begin{align*}
t\mapsto\exp\left\{\eta\int^t_0M_s\dif W_s-\frac{\eta^2}{2}\int^t_0|M_s|^2\dif s\right\}
=:\cE(M)(t)
\end{align*}
is an exponential martingale.
Moreover, by (\ref{ML1}) and Young's inequality
$$
\frac{\eta^2}{2}|M_s|^2-\frac{\eta C_0}{2}X_s^\a+C_1\eta\leq C_{\a,\eta}.
$$
Therefore,
\begin{align*}
\mE \exp\left\{\eta\int^t_0(\frac{C_0}{2}X^\a_s+Y_s)\dif s\right\}
\leq e^{\eta x_0}\mE\left(\cE(M)(t)\cdot\exp\{C_{\a,\eta}t\}\right)= e^{\eta x_0}\cdot\exp\{C_{\a,\eta} t\}.
\end{align*}
The proof is thus complete.
\end{proof}

\subsection{Proof of the Derivative Flow Equation}

In this subsection, we prove (\ref{Der}). Note that (\ref{Lim1})
can be proved similarly.
\bl\label{OO3}
For any $T>0$, there exists a constant $C_{N,T}>0$ such that for each $\om$ and
$\u_0\in\mH^1$
\begin{align*}
\sup_{t\in[0,T]}\|\u(t,\om)\|^2_{\mH^1}+\int^T_0\|\u(t,\om)\|^2_{\mH^2}
\leq C_{N,T}\Big(1+\|\u_0\|^6_{\mH^1}
+\sup_{t\in[0,T]}\|\w(t,\om)\|^{12}_{\mH^5}\Big).
\end{align*}
\el
\begin{proof}
Following the proof of Lemma \ref{LEE4}, let us give different estimates for $J_i, i=1,2,3,4$.

For $J_1$, by (\ref{PP4}) we have
\begin{align*}
J_1\leq -2\|\nabla\v(t)\|^2_{\mH^0}+2\|\Delta\w(t)\|_{\mH^0}\cdot\|\v(t)\|_{\mH^0}.
\end{align*}

For $J_2$, by the Sobolev inequality (\ref{Sob}) and Young's inequality  we have
\begin{align*}
J_2&= -2\<\w(t),((\v(t)+\w(t))\cdot\nabla)(\v(t)+\w(t))\>_{\mH^0}\\
&= 2\<\nabla\w(t),(\v(t)+\w(t))\otimes(\v(t)+\w(t))\>_{\mH^0}\\
&\leq 2\|\nabla\w(t)\|_{L^\infty}\|\v(t)+\w(t)\|^2_{\mH^0}\\
&\leq C\|\nabla\w(t)\|_{L^\infty}\|\v(t)+\w(t)\|^2_{L^4}\\
&\leq \|\v(t)+\w(t)\|^4_{L^4}+C\|\nabla\w(t)\|^2_{L^\infty}.
\end{align*}

For $J_3$, we have
\begin{align*}
J_3\leq-2\|\v(t)+\w(t)\|^4_{L^4}+N\|\v(t)+\w(t)\|^2_{\mH^0}.
\end{align*}

For $J_4$, by (\ref{Sob}) and Young's inequality we have
\begin{align*}
J_3&\leq 2\|\w(t)\|_{L^\infty}\cdot\|\v(t)+\w(t)\|^3_{L^3}\\
&\leq 2\|\w(t)\|_{L^\infty}\cdot\|\v(t)+\w(t)\|^3_{L^4}\\
&\leq -\|\v(t)+\w(t)\|^4_{L^4}+C\|\w(t)\|^4_{L^\infty}.
\end{align*}

Hence
$$
\frac{\dif}{\dif t}\|\v(t)\|^2_{\mH^0}\leq C_N\|\v(t)\|^2_{\mH^0}+C(\|\w(t)\|^4_{\mH^2}
+\|\w(t)\|^2_{\mH^3}).
$$

By Gronwall's inequality, we get
\begin{align}
\sup_{t\in[0,T]}\|\v(t)\|^2_{\mH^0}\leq C_{N,T}\Big(\|\v_0\|^2_{\mH^0}
+\sup_{t\in[0,T]}(\|\w(t)\|^4_{\mH^2}+\|\w(t)\|^2_{\mH^3})\Big).\label{OO2}
\end{align}

Using the similar calculations as in the proof of Lemma \ref{LEE4}, one  finds that
\begin{align*}
\frac{\dif}{\dif t}\|\v(t)\|^2_{\mH^1}\leq-\frac{1}{8}\|\v(t)\|^2_{\mH^2}
+C_N(1+\|\w(t)\|^4_{\mH^5})\cdot(1+\|\v(t)\|^6_{\mH^0}),
\end{align*}
which together with (\ref{OO2}) gives the desired estimate.
\end{proof}
For $\v_0\in\mH^1$, let us consider a small perturbation of
the initial values given by $\u_\epsilon(0)=\u_0+\epsilon\v_0$.
The corresponding solution of Eq. (\ref{NSE00}) is denoted by $\u_\epsilon(t)$.

Set
$$
\v_\eps(t):=(\u_\epsilon(t)-\u(t))/\epsilon.
$$
Then $\v_\eps(t)$ satisfies
\begin{align*}
\v_\eps'(t)&= \Delta\v_\eps(t)-\sP[(\u_\eps(t)\cdot\nabla)\v_\eps(t)]-\sP[(\v_\eps(t)\cdot\nabla)\u(t)]\\
&\quad -\sP[g_N(|\u_\eps(t)|^2)\v_\eps(t)]-\sP[(g_N(|\u_\eps(t)|^2)-g_N(|\u(t)|^2))/\eps\cdot\u(t)],
\end{align*}
 with  initial value $\v_\eps(0)=\v_0$.

We have:
\bl\label{OO4}
For any $T>0$, there is a constant $C_{N,T}>0$ such that for any $\eps\in(0,1)$
$$
\sup_{t\in[0,T]}\|\v_\eps(t)\|^2_{\mH^1}+\int^T_0\|\v_\eps(t)\|^2_{\mH^2}\dif t\leq C_{T,N}.
$$
\el
\begin{proof}
As in the proof of Lemma \ref{Le8}, we have
\begin{align*}
\frac{\dif }{\dif t}\|\v_\eps(t)\|^2_{\mH^1}&\leq -\|\v_\eps(t)\|^2_{\mH^2}
+2\|\v_\eps(t)\|^2_{\mH^1}+C\|\u_\eps(t)\|_{L^6}^2\cdot\|\nabla\v_\eps(t)\|^2_{L^3}\\
&\quad +C\|\v_\eps(t)\|^2_{L^\infty}\cdot\|\u(t)\|^2_{\mH^1}
+C\|\v_\eps(t)\|^2_{L^6}\cdot(\|\u_\eps(t)\|^4_{L^6}+\|\u(t)\|^4_{L^6})\\
&\leq -\frac{1}{2}\|\v_\eps(t)\|^2_{\mH^2}+C(\|\u_\eps(t)\|_{\mH^1}^4+\|\u(t)\|_{\mH^1}^4+1)
\cdot\|\v_\eps(t)\|^2_{\mH^1},
\end{align*}
which together with Lemma \ref{OO3} gives the desired estimate.
\end{proof}

We are now in a position to prove  (\ref{Der2}). Set
$$
\je(t):=\v_\epsilon(t)-\cJ_{0,t}\v_0,
$$
where $\cJ_{0,t}\v_0$ satisfies (\ref{Der}).

By Taylor's formula, we have
\begin{align*}
g_N(|\u_\eps(t)|^2)-g_N(|\u(t)|^2)
&= g'_N(|\u(t)|^2)(|\u_\eps(t)|^2-|\u(t)|^2)
+g''_N(\theta)(|\u_\eps(t)|^2-|\u(t)|^2)^2/2\\
&= \eps^2\cdot g'_N(|\u(t)|^2)|\v_\eps(t)|^2+2g'_N(|\u(t)|^2)\<\v_\eps(t),\u(t)\>_{\mR^3}\\
&\quad +g''_N(\theta)(|\u_\eps(t)|^2-|\u(t)|^2)^2/2,
\end{align*}
where $\theta$ takes some value between $|\u(t)|^2$ and $|\u_\eps(t)|^2$.

Thus, it is not hard to see that $\je(t)$ satisfies
\begin{align*}
\je'(t)=\Delta\je(t)-\sum_{i=1}^8J_i(t),
\end{align*}
where
\begin{align*}
J_1(t)&:= \eps\cdot\sP[(\v_\eps(t)\cdot\nabla)\v_\eps(t)],\\
J_2(t)&:= \sP[(\u(t)\cdot\nabla)\je(t)],\\
J_3(t)&:= \sP[(\je(t)\cdot\nabla)\u(t)],\\
J_4(t)&:= \sP[(g_N(|\u_\eps(t)|^2)-g_N(|\u(t)|^2))\cdot\v_\eps(t)],\\
J_5(t)&:= \sP[g_N(|\u(t)|^2)\cdot\je(t)],\\
J_6(t)&:= \eps\cdot \sP[g'_N(|\u(t)|^2)|\v_\eps(t)|^2\cdot\u(t)],\\
J_7(t)&:= 2\sP[g'_N(|\u(t)|^2)\<\je(t),\u(t)\>_{\mR^3}\cdot\u(t)],\\
J_8(t)&:= \sP[g''_N(\theta)(|\u_\eps(t)|^2-|\u(t)|^2)^2\cdot\u(t)/\eps].
\end{align*}
By the chain rule and Young's inequality, we have
\begin{align*}
\frac{\dif }{\dif t}\|\je(t)\|^2_{\mH^1}\leq-\|\je(t)\|^2_{\mH^2}+2\|\je(t)\|^2_{\mH^0}
+C\sum_{i=1}^8\|J_i(t)\|^2_{\mH^0}.
\end{align*}
Here and below, the constant $C$ is independent of $\eps$.

For $J_1(t)$, we have
$$
\|J_1(t)\|^2_{\mH^0}\leq C \eps^2\cdot\|\v_\eps(t)\|^2_{\mH^1}\cdot\|\v_\eps(t)\|^2_{\mH^2}.
$$

For $J_2(t)$, we have
$$
\|J_2(t)\|^2_{\mH^0}\leq C\|\u(t)\|^2_{\mH^2}\cdot\|\je(t)\|^2_{\mH^1}.
$$

For $J_3(t)$, we have
\begin{align*}
\|J_3(t)\|^2_{\mH^0}&\leq \|\u(t)\|^2_{\mH^1}\cdot\|\je(t)\|^2_{L^\infty}
\leq C\|\u(t)\|^2_{\mH^1}\cdot\|\je(t)\|_{\mH^1}\cdot\|\je(t)\|_{\mH^2}\leq\\
&\leq C\|\u(t)\|^4_{\mH^1}\cdot\|\je(t)\|^2_{\mH^1}+\frac{1}{4}\|\je(t)\|^2_{\mH^2}.
\end{align*}

For $J_4(t)$, we have
\begin{align*}
\|J_4(t)\|^2_{\mH^0}&\leq C\eps^2\cdot\|\v_\eps(t)\|^4_{L^6}\cdot(\|\u_\eps(t)\|^2_{L^6}+
\|\u(t)\|^2_{L^6})\\
&\leq C\eps^2\cdot\|\v_\eps(t)\|^4_{\mH^1}\cdot(\|\u_\eps(t)\|^2_{\mH^1}+
\|\u(t)\|^2_{\mH^1}).
\end{align*}

For $J_5(t)$, we have
\begin{align*}
\|J_5(t)\|^2_{\mH^0}\leq C\|\u(t)\|^4_{L^6}\cdot\|\je(t)\|^2_{L^6}
\leq C\|\u(t)\|^4_{\mH^1}\cdot\|\je(t)\|^2_{\mH^1}.
\end{align*}

For $J_6(t)$, we have
$$
\|J_6(t)\|^2_{\mH^0}\leq C\eps^2\cdot\|\v_\eps(t)\|^2_{\mH^1}\cdot\|\u(t)\|^2_{\mH^1}.
$$

For $J_7(t)$, we have
$$
\|J_7(t)\|^2_{\mH^0}\leq C\cdot\|\je(t)\|^2_{\mH^1}\cdot\|\u(t)\|^4_{\mH^1}.
$$

For $J_8(t)$, we have
\begin{align*}
\|J_8(t)\|^2_{\mH^0}&\leq C\eps^2\cdot\||\v_\eps(t)|^2\cdot(|\u_\eps(t)|^2+|\u(t)|^2)\|^2_{\mH^0}\\
&\leq C\eps^2\cdot\|\v_\eps(t)\|^4_{L^4}\cdot(\|\u_\eps(t)\|^4_{L^4}+\|\u(t)\|^4_{L^4})\\
&\leq C\eps^2\cdot\|\v_\eps(t)\|^4_{\mH^1}\cdot(\|\u_\eps(t)\|^4_{\mH^1}+\|\u(t)\|^4_{\mH^1}).
\end{align*}

Combining the above calculations and Lemmas \ref{OO3} and \ref{OO4} yields that
\begin{align*}
\frac{\dif }{\dif t}\|\je(t)\|^2_{\mH^1}\leq C\eps^2(1+\|\v_\eps(t)\|^2_{\mH^2})+
C(1+\|\u(t)\|^2_{\mH^2})\cdot\|\je(t)\|^2_{\mH^1}.
\end{align*}
By Gronwall's inequality, we get
$$
\|\je(t)\|^2_{\mH^1}\leq C\eps^2\left(1+\int^t_0\|\v_\eps(s)\|^2_{\mH^2}\dif s\right)\cdot
\exp\left\{C+C\int^t_0\|\u(s)\|^2_{\mH^2}\dif s\right\},
$$
which together with Lemmas \ref{OO3} and \ref{OO4} clearly gives  (\ref{Der2}).

\vspace{12mm}

{\bf Acknowledgements:}

\vspace{3mm}

This work was done  while the second named author was a fellow of
Alexander-Humboldt Foundation in Bielefeld University.
He is very grateful to the generous support of Alexander-Humboldt Foundation.
Financial supports of the DFG through SFB 701
and ARC Discovery grant DP0663153 of Australia are also gratefully acknowledged.
The authors would also like to thank the referees for their useful suggestions.

\end{document}